\documentclass[a4paper,11pt,reqno]{amsart}
\newif\iftemplatejapanese
\templatejapanesefalse   % false: 英語専用。LuaTeX-ja と日本語フォントを読み込まない。
\newif\iftemplatelegacy
\templatelegacyfalse     % false: A--Z の短い互換マクロを読み込まない。
\newif\ifdraftnotes
\draftnotesfalse         % false: \note{...} を非表示。

\usepackage[T1]{fontenc}
\usepackage{lmodern}
\usepackage{textcase}

\usepackage{amssymb,mathtools,mathrsfs}
\usepackage[a4paper,margin=30mm,headheight=14pt]{geometry}
\usepackage{microtype}
\usepackage{xcolor}
\usepackage{enumitem}
\setlist[enumerate,1]{label=\textup{(\roman*)},ref=\roman*,leftmargin=*}
\usepackage{comment}
\usepackage{needspace}
\usepackage[unicode,hidelinks]{hyperref}
\usepackage{bookmark}
\usepackage{aliascnt} % MUST precede cleveref; theorem definitions follow cleveref.
\usepackage[nameinlink,capitalize,noabbrev]{cleveref}
\theoremstyle{plain}
\newtheorem{thm}{Theorem}[section]
\newtheorem*{thm*}{Theorem}
\crefname{thm}{Theorem}{Theorems}
\Crefname{thm}{Theorem}{Theorems}

\newcommand{\DeclareSharedTheorem}[3]{%
  \newaliascnt{#1}{thm}%
  \newtheorem{#1}[#1]{#2}%
  \aliascntresetthe{#1}%
  \crefname{#1}{#2}{#3}%
  \Crefname{#1}{#2}{#3}%
}
\DeclareSharedTheorem{lem}{Lemma}{Lemmas}
\DeclareSharedTheorem{prop}{Proposition}{Propositions}
\DeclareSharedTheorem{cor}{Corollary}{Corollaries}
\DeclareSharedTheorem{claim}{Claim}{Claims}
\DeclareSharedTheorem{conj}{Conjecture}{Conjectures}
\newtheorem*{conj*}{Conjecture}

\theoremstyle{definition}
\DeclareSharedTheorem{dfn}{Definition}{Definitions}
\DeclareSharedTheorem{eg}{Example}{Examples}
\DeclareSharedTheorem{prob}{Problem}{Problems}
\DeclareSharedTheorem{setting}{Setting}{Settings}
\DeclareSharedTheorem{question}{Question}{Questions}
\DeclareSharedTheorem{assump}{Assumption}{Assumptions}
\DeclareSharedTheorem{hypo}{Hypothesis}{Hypotheses}
\DeclareSharedTheorem{nota}{Notation}{Notations}
\DeclareSharedTheorem{conv}{Convention}{Conventions}

\newtheorem*{nota-conv}{Notation and Convention}

\theoremstyle{remark}
\DeclareSharedTheorem{rem}{Remark}{Remarks}
\DeclareSharedTheorem{notethm}{Note}{Notes}
\newenvironment{pf}[1][\proofname]{\begin{proof}[#1]}{\end{proof}}

\theoremstyle{plain}

\crefname{athm}{Theorem}{Theorems}
\Crefname{athm}{Theorem}{Theorems}
\newaliascnt{acor}{athm}

\aliascntresetthe{acor}
\crefname{acor}{Corollary}{Corollaries}
\Crefname{acor}{Corollary}{Corollaries}

\newcommand{\ZZ}{\mathbb{Z}}
\newcommand{\QQ}{\mathbb{Q}}
\newcommand{\RR}{\mathbb{R}}
\newcommand{\CC}{\mathbb{C}}
\newcommand{\PP}{\mathbb{P}}
\newcommand{\HH}{\mathbb{H}}

\newcommand{\FF}{\mathbb{F}}

\newcommand{\cC}{\mathcal{C}}

\newcommand{\cH}{\mathcal{H}}

\newcommand{\cO}{\mathcal{O}}
\newcommand{\cP}{\mathcal{P}}
\newcommand{\cQ}{\mathcal{Q}}
\newcommand{\cR}{\mathcal{R}}

\DeclareMathOperator{\Aut}{Aut}

\DeclareMathOperator{\Stab}{Stab}

\DeclareMathOperator{\NS}{NS}
\DeclareMathOperator{\Pic}{Pic}
\DeclareMathOperator{\Amp}{Amp}

\DeclareMathOperator{\rk}{rk}

\DeclareMathOperator{\Isom}{Isom}

\DeclareMathOperator{\MW}{MW}

\DeclareMathOperator{\Nef}{Nef}

\DeclareMathOperator{\vcd}{vcd}

\newcommand{\Orth}{\mathrm O}
\newcommand{\SSph}{\mathbb{S}}
\newcommand{\rt}{\mathrm{rt}}
\newcommand{\qi}{\simeq_{\mathrm{qi}}}
\newcommand{\qs}{\simeq_{\mathrm{qs}}}
\newcommand{\sier}{Sierpi\'nski}
\newcommand{\mobius}{M\"obius}
\DeclareMathOperator{\Cay}{Cay}
\DeclareMathOperator{\Cusp}{Cusp}
\DeclareMathOperator{\diam}{diam}
\DeclareMathOperator{\charac}{char}

\DeclareMathOperator{\spanop}{span}
\DeclareMathOperator{\Km}{Km}

\theoremstyle{plain}
\DeclareSharedTheorem{theorem}{Theorem}{Theorems}
\DeclareSharedTheorem{proposition}{Proposition}{Propositions}
\DeclareSharedTheorem{lemma}{Lemma}{Lemmas}
\DeclareSharedTheorem{corollary}{Corollary}{Corollaries}
\newtheorem*{introTheoremA}{Theorem A}
\newtheorem*{introTheoremB}{Theorem B}
\newtheorem*{introTheoremC}{Theorem C}
\newtheorem*{introTheoremD}{Theorem D}
\theoremstyle{definition}
\DeclareSharedTheorem{definition}{Definition}{Definitions}
\DeclareSharedTheorem{example}{Example}{Examples}
\DeclareSharedTheorem{remark}{Remark}{Remarks}
\DeclareSharedTheorem{problem}{Problem}{Problems}
\DeclareSharedTheorem{introdefinition}{Definition}{Definitions}
\newtheorem*{introquestion}{Question}
\newtheorem*{question*}{Question}

\makeatletter
\AtBeginDocument{%
  \ifx\@pdftitle\@empty
    \hypersetup{pdftitle={\@title}}%
  \fi
  \ifx\@pdfauthor\@empty
    \hypersetup{pdfauthor={\authors}}%
  \fi
}
\makeatother

\numberwithin{equation}{section}
\setlist[enumerate,1]{label=\textup{(\roman*)},ref=\textup{(\roman*)},
  leftmargin=2.2em,itemsep=2pt,topsep=4pt}

\title[Quasi-isometry classes of automorphism groups]{Quasi-isometry classes of automorphism groups of K3 surfaces via Sierpi\'nski spaces}
\author{Kohei Kikuta}
\address{Department of Mathematics, Graduate School of Science,
The University of Osaka, Toyonaka, Osaka 560-0043, Japan}
\email{kikuta@math.sci.osaka-u.ac.jp}
\date{}
\hypersetup{pdflang={en},
  pdfsubject={Quasi-isometry classes of automorphism groups of K3 surfaces}}

\begin{document}
\begin{abstract}
We study quasi-isometry classes of automorphism groups of K3 surfaces using Sierpi\'nski spaces. 
We introduce K3 surfaces of Sierpi\'nski type and prove that their limit sets are round Sierpi\'nski spaces.
For such K3 surfaces, the quasi-isometry class of the automorphism group is equivalent to 
the rational isometry class of the N\'eron--Severi lattice together with its integral roots. 
As an application of this rigidity, we construct infinitely many quasi-isometry classes of hyperbolic automorphism groups in Picard rank 4.
This reveals a new phenomenon in Picard rank 4, sharply contrasting with the finiteness of quasi-isometry classes in Picard ranks at most 3 and at least 6.

\end{abstract}
\maketitle

\setcounter{tocdepth}{1}
%\tableofcontents

% !TEX root = ../main.tex
% Integrated from section1_introduction.tex; full change log: ../CHANGELOG-ja.md.
\section{Introduction}\label{sec:introduction}

\subsection{Automorphism groups of K3 surfaces}
\label{subsec:intro-automorphisms}

The study of automorphism groups of K3 surfaces is a classical subject
in algebraic geometry. These groups range from finite groups to infinite
groups containing nonabelian free subgroups. We first recall several developments for complex K3 surfaces, 
which provide the background for the present paper.

The global Torelli theorem of Piatetski-Shapiro--Shafarevich
\cite{PS71} relates automorphisms of a K3 surface to Hodge isometries
of its second cohomology preserving the ample cone. Together with the
surjectivity of the period map, established by Kulikov and Todorov
\cite{Kul77,Tod80}, it gives a way to pass between geometric questions
and lattice-theoretic constructions. Nikulin developed the theory of
discriminant forms and primitive embeddings \cite{Nik80}, as well as
the study of reflection groups associated with N\'eron--Severi lattices \cite{Nik83}. 
These methods became fundamental tools for studying
automorphism groups. Using hyperbolic reflection groups,
Vinberg explicitly determined the automorphism groups
of two singular K3 surfaces \cite{Vin83}.
In 1985, Sterk proved the cone conjecture: the action on the
rational hull of the nef cone admits a rational polyhedral fundamental
domain \cite{Ste85}. In particular, the automorphism group is finitely
generated. 
The corresponding result in characteristic
$p>2$ was later established by Lieblich--Maulik
\cite{LM18}.
In 1988, Mukai classified finite groups acting symplectically
on complex K3 surfaces: they are precisely the subgroups of eleven
maximal groups \cite{Muk88}. More recently, Brandhorst--Hofmann
classified finite group actions up to deformation \cite{BH23}.

We next recall the distinction between elementary
and non-elementary automorphism groups.
%We next distinguish automorphism groups according to whether their natural action on hyperbolic space is elementary or non-elementary.
%There is another useful distinction, concerning the full automorphism group. 
Let $\NS(X)$ denote the N\'eron--Severi lattice of $X$,
and put $\rho(X):=\rk\NS(X)$. 
The intersection form has signature $(1,\rho(X)-1)$, and the positive cone determines
a real hyperbolic space $\HH_X$ of dimension $\rho(X)-1$. 
%For $\rho(X)\geq2$, 
The natural representation
\[
 \varrho_X:\Aut(X)\longrightarrow\Orth^+(\NS(X))
 \subset\Isom(\HH_X),%\qquad
 %g\longmapsto(g^{-1})^*|_{\NS(X)},
\]
has finite kernel. 
%Write $\Gamma_X$ for its image. 
We call %this action, and hence 
$\Aut(X)$ \emph{elementary} if $\Aut(X)$ has a
finite orbit in the ideal compactification $\HH_X\cup\partial_\infty\HH_X$. Equivalently,
its limit set
\[
 \Lambda_X:=\overline{\Aut(X)o}\cap\partial_\infty\HH_X,
 \qquad o\in\HH_X,
\]
has at most two points. Here the closure is taken in the ideal
compactification. For these groups, elementary is equivalent to
virtually abelian. It includes finite groups, groups preserving a
parabolic fixed point, and virtually cyclic groups preserving a pair of boundary points. 
%When $\rho(X)=1$, the automorphism group is finite.

The elementary case is now well understood in terms of the N\'eron--Severi lattice. 
The finite case was classified
at the level of N\'eron--Severi lattices by Nikulin and
Vinberg \cite{Nik83,Nik84,Vin07}, and the corresponding
automorphism groups were further studied by Kond\=o
\cite{Kon89}; see also \cite{BKM26} and the references
therein. Oguiso \cite{Ogu07} showed that 
every zero-entropy subgroup is virtually abelian and that an infinite zero-entropy subgroup preserves an elliptic fibration. 
%every null-entropy subgroup is almost abelian and that an infinite null-entropy subgroup preserves an elliptic fibration. 
Brandhorst--Mezzedimi \cite{BM23}, and
independently Yu \cite{Yu25} for K3 surfaces, classified
the N\'eron--Severi lattices in the zero-entropy case with infinite automorphism group.
Brandhorst--Kirschmer--Mezzedimi subsequently completed
the classification by treating elementary groups of
hyperbolic type \cite{BKM26}. Thus, in Picard rank at least
3, the N\'eron--Severi lattices giving elementary automorphism groups form an explicit finite list. 
%This classification concerns the full group, whereas Mukai's theorem concerns its finite symplectic subgroups.

For non-elementary automorphism groups, several complementary
approaches have been developed. The dynamical approach studies entropy,
invariant measures, and orbits on the surface. Cantat's work on K3
automorphisms \cite{Can01} and Oguiso's Tits alternative for projective
hyperk\"ahler manifolds \cite{Ogu06} are fundamental contributions in
this direction. In particular, the latter gives a dichotomy between
virtually abelian groups and groups containing nonabelian free
subgroups. A different approach seeks explicit generators through
lattice embeddings and chamber decompositions. Borcherds' method
\cite{Bor98}, developed algorithmically by Shimada \cite{Shi15},
has been effective for many individual K3 surfaces.

Hyperbolic geometry provides a framework for studying the non-elementary case collectively. 
Totaro %systematically 
developed the relation between
cone conjectures and hyperbolic geometry \cite{Tot10,Tot12}, including
examples whose automorphism groups are not commensurable with
arithmetic groups. Kurnosov--Yasinsky studied automorphism groups of
projective hyperk\"ahler manifolds using actions on CAT(0) spaces
\cite{KY23}. More recently, the author \cite{Kik26} and, independently
over $\CC$, Takatsu \cite{Tak25} proved geometrical finiteness of the action on $\HH_X$
Consequently, a non-elementary automorphism group is relatively hyperbolic with respect to its infinite genus-one-fibration stabilizers. %, up to conjugacy. 
Such a stabilizer is virtually abelian, with rank equal to the Mordell--Weil rank of the Jacobian fibration. In particular,
\[
 \Aut(X)\text{ is hyperbolic}
 \quad\Longleftrightarrow\quad
 \rk\MW(\pi)\leq1
 \text{ for every genus-one fibration }\pi.
\]
%This condition is vacuous when there is no genus-one fibration.
See %\cite[Theorem~1.3]{FOY26} and 
Theorem~\ref{thm:k3-hyperbolic} for details. %for the criterion and its extension to arbitrary characteristic.

In this approach, the limit set $\Lambda_X$ is an important invariant of $\Aut(X)$. %the hyperbolic representation. 
Its topology and metric geometry relate the structure of the group to the boundary of the nef cone.
The limit set is naturally identified with the Bowditch boundary of $\Aut(X)$ with respect to the full parabolic peripheral family. 
%For a convex-cocompact action, it identifies with the Gromov boundary of the group; 
%in the geometrically finite case, it is instead related to the Bowditch boundary for the full parabolic family. 
%This distinction is relevant when a fibration stabilizer is virtually cyclic. 
The present paper uses limit sets to study automorphism groups from the viewpoint of quasi-isometry.

\subsection{Quasi-isometry classes of automorphism groups}
\label{subsec:intro-qi}

The results above give structural theorems and many explicit examples.
A general classification of non-elementary automorphism groups,
however, remains largely open. In this paper, we formulate and study
their classification up to quasi-isometry. This provides a coarse
classification problem in which one can compare groups arising from
different K3 surfaces without first determining their generators and
relations.

Quasi-isometry is a fundamental equivalence relation in geometric
group theory. A finitely generated group is equipped with a word
metric, whose quasi-isometry class is independent of the choice of
a finite generating set. Roughly speaking, a quasi-isometry compares
distances up to a multiplicative constant and a bounded additive error.
Passage to a finite-index subgroup or to a quotient by a finite normal
subgroup does not change this class. Although coarse, this equivalence
preserves substantial information, including growth type, the number
of ends, and hyperbolicity. Among hyperbolic groups, it also preserves
elementarity and the homeomorphism type of the Gromov boundary. The
boundary construction of canonical splittings over two-ended subgroups
gives a further manifestation of quasi-isometry invariance for
one-ended hyperbolic groups \cite{Bow98}. 
%which are not virtually Fuchsian \cite{Bow98}. 
%We refer to \cite{BH99} for general background.

It is natural to begin with hyperbolic automorphism groups, since
hyperbolicity itself is invariant under quasi-isometry. 
From now on, unless otherwise stated, all K3 surfaces are
over algebraically closed fields of arbitrary characteristic.
For $\rho\in\{1,\ldots,20,22\}$, set
\begin{equation}\label{eq:intro-Q}
 \cQ_\rho:=
 \left\{
 [\Aut(X)]_{\mathrm{qi}}\ \middle|\
 \begin{array}{l}
 \rho(X)=\rho,\\
 \Aut(X)\text{ is non-elementary hyperbolic}
 \end{array}
 \right\},
 \qquad Q_\rho:=\#\cQ_\rho.
\end{equation}
Here the K3 surface and its algebraically closed ground field both vary, 
and $[\Aut(X)]_{\mathrm{qi}}$ denotes the quasi-isometry class of $\Aut(X)$. 
For hyperbolic groups, \emph{non-elementary} means not virtually
cyclic; this agrees with non-elementarity of $\Aut(X)$ defined the above. 
%the natural action in the present setting. 
In particular, finite groups are excluded.
For $\rho\leq20$, the set $\cQ_\rho$ is unchanged if we restrict to
complex K3 surfaces. %, as explained in Subsection~\ref{subsec:hyperbolic-qi-classes}.

The low-rank cases and the supersingular case give
\[
 Q_1=Q_2=Q_{22}=0,\qquad Q_3=2.
\]
The two classes in Picard rank 3 are represented by a nonabelian
free group of finite rank and by the fundamental group of a closed
hyperbolic surface. 
\begin{comment}
Both occur as quasi-isometry classes of K3 automorphism groups: 
a very general Wehler surface gives the first, and a K3 surface with root-free anisotropic N\'eron--Severi lattice
\[
 \langle4\rangle\oplus\langle-8\rangle\oplus\langle-12\rangle
\]
gives the second. 
\end{comment}
Details are given in Subsection~\ref{subsec:hyperbolic-qi-classes}.
%Details and explicit constructions are given in Subsection~\ref{subsec:hyperbolic-qi-classes}.

In higher Picard rank, Fujiwara--Oguiso--Yu \cite[Theorem~1.8]{FOY26}
proved that exactly $166$ N\'eron--Severi lattices occur for complex
K3 surfaces of Picard rank at least 6 with non-elementary hyperbolic automorphism group. 
Isometric N\'eron--Severi lattices give quasi-isometric automorphism groups. 
Their theorem therefore implies
\[
 \sum_{\rho=6}^{20}Q_\rho\leq166.
\]
Different lattices may contribute the same quasi-isometry class, so
this bound does not determine the individual numbers $Q_\rho$.
Our first main result concerns Picard rank 4.
%Our main infinitude result concerns Picard rank 4.

\begin{introTheoremA}[Theorem~\ref{thm:Q-four-infinite}]
We have $Q_4=\infty$. More precisely, there exist infinitely many
complex K3 surfaces of \sier\ type and Picard rank 4 whose
automorphism groups are non-elementary hyperbolic and pairwise
non-quasi-isometric.
\end{introTheoremA}

Thus infinitely many classes remain even after passing to this coarse
equivalence relation. Moreover, all groups in the construction have
Gromov boundary homeomorphic to a \sier\ carpet. Boundary topology
alone does not distinguish them; the quasisymmetric rigidity discussed below is the key to the proof. 
%The present paper does not determine whether $Q_5$ is finite.

\subsection{\texorpdfstring{K3 surfaces of \sier\ type}{K3 surfaces of Sierpinski type}}
\label{subsec:intro-sierpinski}

\sier\ carpets are basic objects in the study of boundaries of
hyperbolic groups. 
Kapovich--Kleiner \cite[Theorem~4]{KK00} showed that
the one-dimensional boundary of a non-elementary hyperbolic group
admitting no nontrivial splitting over a finite or virtually cyclic
subgroup is homeomorphic to a circle, a \sier\ carpet, or the Menger curve.
Among these possibilities, the \sier\ carpet is the only planar case other than the circle.
This makes it a natural boundary to consider in the quasi-isometry classification.

Our interest also comes from a metric rigidity phenomenon. 
A $m$-dimensional round \sier\ space is obtained from $\SSph^{m+1}$ by removing a countably infinite dense union of round open balls with pairwise disjoint closures and diameters tending to zero. 
The case $m=1$ is a round \sier\ carpet. 
Bonk--Kleiner--Merenkov \cite[Theorem 1.1]{BKM09} proved that a quasisymmetric map from a round \sier\ space of spherical measure zero onto another round \sier\ space in the same sphere is the restriction of a \mobius\ transformation. 
\begin{comment}
These sets are Schottky sets, and
Bonk--Kleiner--Merenkov \cite[Theorem 1.1]{BKM09} proved that a quasisymmetric map from a Schottky set of spherical measure zero onto another Schottky set in the same sphere is the restriction of a \mobius\ transformation, provided the sphere has dimension at least two.
\end{comment}
Thus rigidity depends on the round embedding and the measure-zero
condition, not only on the underlying topology.

We introduce a class of K3 surfaces for which this geometry is
encoded by intersection numbers of $(-2)$-curves.
%We use the notation $\NS(X)$, $\HH_X$, and $\Lambda_X$ from Subsection~\ref{subsec:intro-automorphisms} also in arbitrary characteristic. 
%Here a \emph{$(-2)$-curve} means a smooth rational curve on $X$, and $\cR_X$ denotes the set of such curves.
We denote by $\cR_X$ the set of $(-2)$-curves on $X$.
 
\begin{definition}\label{def:intro-sierpinski}
A K3 surface $X$ is \emph{of \sier\ type} 
if the following conditions hold:
\begin{enumerate}
\item $\rho(X)\geq4$; 
\item $X$ contains a $(-2)$-curve;
\item $(C,D)\geq3$ for any distinct $C,D\in\cR_X$.
\end{enumerate}
\end{definition}

%The Picard-rank restriction places the corresponding boundary in the range of positive-dimensional \sier\ spaces. 
The definition is algebro-geometric, while its relation
to the limit set is given by the following theorem.

\begin{introTheoremB}[Theorem~\ref{thm:round-limit-set}]
If $X$ is of \sier\ type, then $\Lambda_X$ is a
$(\rho(X)-3)$-dimensional round \sier\ space.
Each complementary component is precisely the
round open ball $B_C~(C\in\cR_X)$, 
and the assignment $C\mapsto B_C$ is a bijection.
%\[
% B_C=\{[v]\in\partial_\infty\HH_X\mid(v,C)<0\},
% \qquad C\in\cR_X.
%\]
\end{introTheoremB}

\begin{comment}
We do not establish the converse. We formulate it as follows in
Subsection~\ref{subsec:round-limit-sets}.
\begin{introquestion}
Let $X$ be a K3 surface with $\rho(X)\geq4$. If its geometric limit
set is a $(\rho(X)-3)$-dimensional round \sier\ space in
$\partial_\infty\HH_X$, must $X$ be of \sier\ type?
\end{introquestion}
\end{comment}

These surfaces satisfy several further algebro-geometric restrictions.

%The geometry of these surfaces is constrained in several other ways.
%For a genus-one fibration $\pi$, write $\MW(\pi)$ for the Mordell--Weil group of its Jacobian fibration.

\begin{introTheoremC}[Theorem~\ref{thm:fundamental}]
Let $X$ be a K3 surface of \sier\ type. Then the following hold.
\begin{enumerate}
\item
$\rho(X)\in\{4,5\}$.
\item
Every genus-one fibration
$\pi:X\to\PP^1$ is elliptic and has no reducible fibers, i.e. its singular fibers are of type $I_1$ or $I\!I$, and
\[
 \MW(\pi)\cong\ZZ^{\rho(X)-2}.
\]
%Here $\MW(\pi)$ denotes the Mordell--Weil group of the Jacobian fibration, as in Section~\ref{sec:preliminaries}.
\item 
$\Lambda_X=\Lambda_X^{\mathrm{nef}}$,
where $\Lambda_X^{\mathrm{nef}}$ denotes the boundary at infinity
of the projectivized nef cone of $X$.
\item
$\Aut(X)$ is non-elementary.
\item
$X$ contains infinitely many $(-2)$-curves.
\item
The following conditions are equivalent:
$\Aut(X)$ is hyperbolic; $\Aut(X)$ is convex-cocompact;
$X$ admits no genus-one fibration; and $\NS(X)$ is anisotropic.
If these conditions hold, then $\rho(X)=4$.
\end{enumerate}
\end{introTheoremC}

To state the quasi-isometry rigidity theorem, let
$
 \Delta_X:=\{r\in \NS(X)\mid r^2=-2\}.
$
A rational isometry $\phi:\NS(X_1)_\QQ\to \NS(X_2)_\QQ$ is \emph{root-preserving} if $\phi(\Delta_{X_1})=\Delta_{X_2}$.
%This is a condition on the integral root sets of the specified lattices; it does not follow from rational isometry alone.

\begin{introTheoremD}[Theorem~\ref{thm:rigidity-sierpinski}]
Let $X_1$ and $X_2$ be K3 surfaces of \sier\ type, possibly over
different algebraically closed fields. 
The following conditions are equivalent:
\begin{enumerate}
\item
$\Aut(X_1)\qi\Aut(X_2)$;
\item
there is a root-preserving rational isometry
$\phi:\NS(X_1)_\QQ\xrightarrow{\sim}\NS(X_2)_\QQ$;
\item
$\rho(X_1)=\rho(X_2)$, and their limit sets are \mobius\ equivalent in the corresponding boundary sphere.
\end{enumerate}
\end{introTheoremD}

The key ingredient in the proof is the quasisymmetric rigidity theorem of Bonk--Kleiner--Merenkov mentioned above.
To prove Theorem~A, 
we apply Theorem~D to an explicit infinite family of
pairwise rationally non-isometric even hyperbolic lattices of rank 4, realized as the Néron--Severi lattices of complex K3 surfaces of \sier\ type. 

\begin{comment}
The proof passes from a quasi-isometry to a quasisymmetry of the
limit sets. Theorem~C ensures that every member of the full parabolic
peripheral family has rank at least two, so this passage also applies
in the relatively hyperbolic case. Since these limit sets are proper,
geometrical finiteness gives the required measure-zero property. The
Bonk--Kleiner--Merenkov theorem then supplies a \mobius\
transformation. By Theorem~B, it permutes the complementary balls
coming from $(-2)$-curves. Their normal vectors have the fixed square
$-2$ and span the rational N\'eron--Severi space, allowing us to
recover the root-preserving rational isometry. The reverse implication
follows from the comparison of ample-cone stabilizers in
Proposition~\ref{prop:root-qi}.
\end{comment}

\begin{comment}
Finally, we apply Theorem~D to an explicit family of rank-four
lattices with Gram matrices
\[
 G_t=
 \begin{pmatrix}
 2&t&t&t\\
 t&-2&4&4\\
 t&4&-2&4\\
 t&4&4&-2
 \end{pmatrix},
 \qquad t>0,\quad t\equiv3\pmod{12}.
\]
These lattices are anisotropic and occur as N\'eron--Severi lattices
of complex K3 surfaces of \sier\ type. Their determinants are
$-108(t^2-4)$. Choosing infinitely many parameters for which
$t^2-4$ have distinct square classes in $\QQ^\times$ excludes rational
isometries and, by Theorem~D, quasi-isometries of the automorphism
groups. This proves Theorem~A without using the converse to Theorem~B.
\end{comment}

\subsection{Organization of the paper}\label{subsec:intro-organization}

Section~\ref{sec:preliminaries} recalls the required facts on K3
automorphism groups, quasi-isometries, relative hyperbolicity, and
\sier\ spaces. Section~\ref{sec:sierpinski-type} introduces K3
surfaces of \sier\ type and proves Theorems~C, B, and~D, in this
order. 
%It also discusses the recognition question following Theorem~B. 
Section~\ref{sec:rank-four} studies the numbers $Q_\rho$
and the finiteness result in higher Picard rank. The final subsection constructs the rank-four family and proves Theorem~A.

\subsection*{AI disclosure}

This paper is part of the author's project on automorphism groups
using geometric group theory. Its starting point was the finiteness
result of Fujiwara--Oguiso--Yu, announced in July 2025, which suggested
that a quasi-isometry classification might be accessible, at least
for hyperbolic automorphism groups. The idea of using \sier\ spaces
is due to the author. The formulations of Theorems~A--D and most of
their proofs are due to the author, while some technical details were
refined through discussions with OpenAI's ChatGPT, including
GPT-6 Astra.
During the preparation of this paper, these tools also assisted with
the search for explicit lattice examples and related computations,
the examination of mathematical arguments, literature searches,
drafting and revising the text and LaTeX code, language editing, and
proofreading.
% Author confirmation required before submission: the following sentence
% is retained from the supplied introduction and is not a verification
% certified by the integration/build checks.
All mathematical statements, proofs, and references in
the final manuscript were independently verified by the author, who
takes full responsibility for the content of this paper.

\subsection*{Acknowledgements}

The author is indebted to Arend Bayer and the University of Edinburgh
for their kind hospitality. This work was supported by a JSPS Overseas
Research Fellowship and JSPS KAKENHI Grant Number 21K13780.

% !TEX root = ../main.tex
% Integrated from section2_preliminaries.tex; full change log: ../CHANGELOG-ja.md.
\section{Preliminaries}\label{sec:preliminaries}

In this section, we collect basic facts on automorphism groups of K3
surfaces, quasi-isometries, relatively hyperbolic groups, and \sier\ spaces.
%Unless otherwise stated, a K3 surface is defined over an algebraically closed field, and all groups considered as metric spaces are finitely generated and endowed with word metrics.

\subsection{Automorphism groups of K3 surfaces}\label{subsec:k3}

We first fix notation for K3 surfaces and their associated lattices.
We refer to \cite{Huy16, CDL25, Dolgachev-Kondo_EnriquesII, Kik26}. 
%for general background.
In this paper, 
a \emph{K3 surface} over an algebraically closed field $k$ of arbitrary characteristic 
is a smooth projective surface $X$ such that
$\omega_X\cong\cO_X$ and $H^1(X,\cO_X)=0$. 
The group of $k$-automorphisms of $X$ is denoted by $\Aut(X)$.
The N\'eron--Severi group
$
 \NS(X):=\Pic(X) %:=\Pic(X)/\Pic^0(X)
$
is a finitely generated free abelian group, and its rank $\rho(X):=\rk\NS(X)$ is called the \emph{Picard rank}. 
The intersection pairing makes it an even
lattice of signature $(1,\rho(X)-1)$, 
i.e. $\NS(X)$ is a hyperbolic lattice. 
For $\mathbb F=\QQ$ or $\RR$, we write
$\NS(X)_\FF:=\NS(X)\otimes_\ZZ\mathbb F$ and use the same notation
$(\ ,\ )$ for the extended pairing. We write $v^2:=(v,v)$.

Let $\cC_X$ be the connected component of
$\{v\in \NS(X)_\RR\mid v^2>0\}$ containing ample classes. 
The \emph{ample cone} $\Amp(X)$ is the cone generated by ample divisor classes, and
the \emph{nef cone} $\Nef(X):=\overline{\Amp(X)}$ is its closure in $\NS(X)_\RR$.
We set
\[
 \Orth^+(\NS(X)):=\{g\in\Orth(\NS(X))\mid g(\cC_X)=\cC_X\},
 \]
 \[
 A_X:=\{g\in\Orth^+(\NS(X))\mid g(\Nef(X))=\Nef(X)\}.
 \]
We use push-forward to obtain a left action:
\begin{equation}\label{eq:representation}
 \varrho_X:\Aut(X)\longrightarrow A_X,\qquad
 f\longmapsto f_*=(f^{-1})^*|_{\NS(X)}.
\end{equation}
Its image is denoted by $\Gamma_X$.

A vector $r\in \NS(X)$ with $r^2=-2$ is called a \emph{root}. We denote
the set of roots by $\Delta_X$. The reflection associated with a root is
$
 s_r(v):=v+(v,r)r
$
for $v\in\NS(X)$. 
The subgroup $W_X\subset\Orth^+(\NS(X))$ generated by these reflections is
called the \emph{Weyl group}. The root hyperplanes divide $\cC_X$ into chambers,
and $W_X$ acts simply transitively on the set of chambers. 
The ample cone is one of these chambers; more precisely, we have,
\begin{equation}\label{eq:ample}
\Amp(X)=\{x\in\cC_X\mid (x,[C])>0
       \text{ for every $(-2)$-curve }C\subset X\}.
\end{equation}
%Indeed, Riemann--Roch implies that one of $r$ and $-r$ is effective for every root $r$. 
%A root need not itself be the class of an irreducible curve. 
We also have
\[
 \Orth^+(\NS(X))=W_X\rtimes A_X.
\]

\begin{comment}
We will distinguish the field-independent notation above from notation
for integral cohomology. If $k=\CC$, then
$H^2(X,\ZZ)$, with its cup-product pairing, is the K3 lattice of signature
$(3,19)$. We write $\Orth^+(H^2(X,\ZZ))$ for the subgroup preserving an
orientation of maximal positive definite three-dimensional subspaces.
The complex structure and a K\"ahler class determine this orientation:
one can use the ordered basis
$\Re\omega,\Im\omega,h$, where
$0\ne\omega\in H^{2,0}(X)$ and $h$ is a K\"ahler class.
The action of $\Aut(X)$ on $H^2(X,\ZZ)$ preserves this orientation and
the Hodge structure. In arbitrary characteristic, we use
$\Orth^+(\NS(X))$ as above; no choice of a Betti cohomology lattice is
implicit in this notation.
\end{comment}

We record explicitly the characteristic hypothesis needed
to compare $\Aut(X)$ with $A_X$.
In positive characteristic, a K3 surface with $\rho(X)=22$ is called \emph{supersingular}. 
Its discriminant is $-p^{2\sigma(X)}$, where $p:=\charac k$ and $1\leq\sigma(X)\leq10$ is its \emph{Artin invariant}. 
Consider the condition
\begin{equation}\label{eq:char}
 \text{$X$ is not supersingular in characteristic $2$.}
\end{equation}

\begin{comment}
\begin{equation}
 \text{$X$ is not supersingular of Artin invariant $10$
 in characteristic $2$.}
\end{equation}
\end{comment}

\begin{theorem}[{\cite{Dolgachev-Kondo_EnriquesII, Ste85, LM18}}]
\label{thm:cone-comparison}
The homomorphism \eqref{eq:representation} has finite kernel.
If $X$ satisfies \eqref{eq:char}, then $\Gamma_X$ has finite index in
$A_X$, and $\Aut(X)$ is finitely generated.
\end{theorem}

The finite-kernel assertion holds in every characteristic; see, for example, \cite[Proposition~8.2.1]{Dolgachev-Kondo_EnriquesII}.
The comparison with the chamber stabilizer $A_X$ follows from the (crystalline) Torelli theorem and cone conjecture; see \cite{Ste85}, \cite[Theorem~6.1 and Section~6]{LM18}.
%, and \cite[Theorem~2.2]{Bra18}. 
%For supersingular surfaces in characteristic $2$, the range $\sigma(X)<10$ is explained in \cite[Remark~1.3]{Bra18}. 
We retain \eqref{eq:char} whenever this comparison is used. 
%The criterion for hyperbolicity in Theorem~\ref{thm:k3-hyperbolic} below will not require it.

A \emph{genus-one fibration} on $X$ is 
a surjective morphism $\pi:X\to\PP^1$ such that $\pi_*\mathcal O_X=\mathcal O_{\mathbb P^1}$. 
%By \cite[Theorem 4.1.1]{CDL25} and the adjunction formula, 
Its generic fiber is a geometrically integral regular curve of arithmetic genus one. 
Such a fibration is called \emph{elliptic} if its generic fiber is smooth, and \emph{quasi-elliptic} otherwise.
%A \emph{genus-one fibration} on $X$ is a surjective morphism $\pi:X\to\PP^1$ with connected fibers whose generic fiber is a geometrically integral curve of arithmetic genus one. 
%It is \emph{elliptic} if the generic fiber is smooth, and \emph{quasi-elliptic} otherwise. 
The latter case occurs only for supersingular K3 surfaces
in characteristics $2$ and $3$.
We write $\MW(\pi)$ for the group of sections of its Jacobian fibration, and call the \emph{Mordell--Weil group} of $\pi$. 
It is a finitely generated abelian group and acts on $X$ by translations along the fibres.
Set the stabilizer of $\pi$
\[
 \Aut(X,\pi):=\{g\in\Aut(X)\mid
 \pi\circ g=\bar g\circ\pi\text{ for some }\bar g\in\Aut(\PP^1)\}.
\]
Then, for the fiber class $e$ of $\pi$, we have
$\Aut(X,\pi)=\{g\in\Aut(X)\mid g^*e=e\}$. 
The Mordell--Weil group $\MW(\pi)$ is a finite-index subgroup of $\Aut(X,\pi)$ in every characteristic; see, for example,  \cite[Proposition~4.4]{KTT26}.

%Let $e\in\NS(X)$ be the fiber class of $\pi$. 
The fiber class of $\pi$ is primitive, nef, and isotropic. 
Conversely, every nonzero primitive nef isotropic
class $e\in\NS(X)$ defines a genus-one fibration $\pi_e:X\to\PP^1$. 
The lattice
\[
 M_e:=e^\perp_{\NS(X)}/\ZZ e
\]
is even negative definite of rank $\rho(X)-2$.
For a negative definite lattice $M$, let $M_{\rt}$ denote the sublattice
generated by its roots. The Shioda--Tate formula gives
\begin{equation}\label{eq:mw-rank}
 \rk\MW(\pi)=\rho(X)-2-\rk\bigl(M_e)_{\rt}\bigr.
\end{equation}
%This also applies to a fibration without a section, with the Mordell--Weil group understood as above; see \cite{CDL25} and \cite[proof of Proposition~4.4]{KTT26}.
Here the roots in $M_e$ span the root lattice of the reducible fibers. %: by the Riemann--Roch theorem, a root orthogonal to the nef fiber class is effective up to sign and is supported in fibers.
For a quasi-elliptic fibration, its Mordell--Weil rank is zero.  %\cite[Proposition~4.5]{KTT26}.

\vspace{4mm}

For $\rho(X)\geq2$, let
$
\HH_X:=\cC_X/\RR_{>0}.
$
We identify $\HH_X$ with the \emph{hyperboloid}
\[
\{x\in\cC_X\mid x^2=1\}
\]
via $[x]\mapsto x/\sqrt{x^2}$.
Equipped with the metric
%\begin{equation}\label{eq:hyperboloid}
$\cosh d_{\HH_X}(x,y)=(x,y)$,
%\end{equation}
it is a model of real hyperbolic $(\rho(X)-1)$-space.
Its \emph{ideal boundary} is
\[
 \partial_\infty\HH_X:=
 (\partial\cC_X\setminus\{0\})/\RR_{>0}.
\]
The action \eqref{eq:representation} gives %a proper
an isometric action on $\HH_X$ with finite kernel. 
%For $\rho(X)=1$, the same construction gives a point and $\Aut(X)$ is finite.
%For $\rho(X)\geq2$, 

\begin{dfn}
\begin{enumerate}
\item
The \emph{limit set} of $\Aut(X)$ is defined by
\[
\Lambda_X:=
\Lambda_{\Aut(X)}:=
\overline{\Aut(X) o}\cap\partial_\infty\HH_X,
%=\overline{\Gamma_X o}\cap\partial_\infty\HH_X,
\]
where the closure is taken in the ideal compactification $\HH_X\cup\partial_\infty\HH_X$ of $\HH_X$, 
and $\Lambda_X$ is independent of the choice of $o\in\HH_X$.
\item 
The \emph{nef boundary} of $X$ is defined by
\[
\Lambda_X^{\mathrm{nef}}
:=\{[v]\in\partial_\infty\HH_X
\mid v\in\Nef(X)\}. 
\]
%, where $[v]$ denotes the positive ray through $v$.
\end{enumerate}
\end{dfn}
Let $\pi:X\to\PP^1$ be a genus-one fibration with
fiber class $e\in\NS(X)$.
Every infinite-order element of $\MW(\pi)$ acts on
$\HH_X$ as a parabolic isometry whose unique fixed
point on $\partial_\infty\HH_X$ is $[e]$.
In particular, if $\rk\MW(\pi)>0$, then $[e]$ is a
parabolic fixed point of $\Aut(X)$ and belongs to
$\Lambda_X\cap\Lambda_X^{\mathrm{nef}}$.

We call $\Aut(X)$ \emph{elementary} if $\Lambda_X$ is
finite, equivalently, if $\#\Lambda_X\leq2$.
This is also equivalent to $\Aut(X)$ having a finite
orbit in the ideal compactification of $\HH_X$,
or to $\Aut(X)$ being virtually abelian of finite rank.
%\cite[Definition~2.1 and Proposition~2.2]{Kik26}.
%elementary Klein群は自動的に有限生成，よって標数2でも大丈夫(そもそもsupersingularなら常にnon-elementaryだから初めから問題はないんだけど)
Moreover, $\Aut(X)$ is non-elementary 
whenever $X$ is supersingular \cite[Corollary~5.4]{BM23}.

We generalize the Baragar--Takatsu theorem to arbitrary
characteristic. 
The result does not require the assumption \eqref{eq:char}.
\begin{thm}
[{Baragar--Takatsu \cite[Theorem 3.3]{BaragarAmpleConeOrbits}, \cite[Theorem~5.28]{Tak25}}]
\label{limit-set-nef-boundary}
Let $X$ be a K3 surface %over an algebraically closed field of arbitrary characteristic, 
with $\rho(X)\geq3$. 
Then
\[
 \Lambda_X^{\mathrm{nef}}=\Lambda_X\sqcup E_X,
\]
where
$
 E_X:=
 \bigl\{[e]\in\partial_\infty\HH_X\mid e\in\NS(X)\text{ is primitive, nef and isotropic,}
                  \rk\MW(\pi_e)=0\bigr\},
$
and every point of $E_X$ is isolated in $\Lambda_X^{\mathrm{nef}}$.
\end{thm}

\begin{pf}
Since $\Aut(X)$ preserves $\Nef(X)$,
we have $\Lambda_X\subseteq\Lambda_X^{\mathrm{nef}}$.
Put $L:=\NS(X)$, $G:=\Orth^+(\NS(X))$, and
$D:=\Nef(X)\cap\HH_X$.
For a primitive nef isotropic class $e\in L$, set
\[
 M_e:=e^\perp_L/\ZZ e,\qquad
 R_e:=(M_e)_{\mathrm{rt}},\qquad
 r(e):=\rho(X)-2-\rk R_e.
\]
The Shioda--Tate formula \eqref{eq:mw-rank} gives
$r(e)=\rk\MW(\pi_e)$.

\noindent(i)
Let $\xi=[e]\in\Lambda_X^{\mathrm{nef}}$ be rational with $r(e)>0$.
The class $e\in\NS(X)$ defines a genus-one fibration $\pi_e$.
Since $\rk\MW(\pi_e)=r(e)>0$, the point $[e]$ is
a parabolic fixed point of $\Aut(X)$, as recalled above.
Hence $[e]\in\Lambda_X$.

\noindent(ii)
Let $\xi=[e]\in\Lambda_X^{\mathrm{nef}}$ be rational with $r(e)=0$.
Choose roots $r_1,\ldots,r_{\rho(X)-2}\in e^\perp_L$
whose images form a basis of $(M_e)_\RR$.
For an ample class $h$ and a sufficiently large integer $N$,
both $Ne+r_i$ and $Ne-r_i$ have positive intersection with $h$.
They have square $-2$, so the Riemann--Roch theorem implies that
they are effective.
Every $[x]\in\Lambda_X^{\mathrm{nef}}\setminus\{[e]\}$
has a unique representative with $(x,e)=1$, and nefness gives
\begin{equation}\label{eq:isolated}
 |(x,r_i)|\leq N\qquad(1\leq i\leq\rho(X)-2).
\end{equation}
To interpret these bounds, choose $f\in L_\RR$ with $f^2=0$
and $(e,f)=1$, and put $V:=(\RR e+\RR f)^\perp\simeq (M_e)_\RR$.
The standard affine chart on
$\partial_\infty\HH_X\setminus\{[e]\}$ is
\[
 V\longrightarrow\partial_\infty\HH_X\setminus\{[e]\},
 \qquad
 v\longmapsto\left[f+v-\frac{v^2}{2}e\right].
\]
Since the images of the $r_i$ form a basis of $(M_e)_\RR$,
the inequalities \eqref{eq:isolated} show that
the preimage of $\Lambda_X^{\mathrm{nef}}\setminus\{[e]\}$
in $V$ is bounded.
Thus $\Lambda_X^{\mathrm{nef}}\setminus\{[e]\}$ is contained
in a compact subset of $\partial_\infty\HH_X\setminus\{[e]\}$,
and hence $[e]$ is isolated in $\Lambda_X^{\mathrm{nef}}$.

If $\Aut(X)$ is non-elementary, its limit set is perfect,
so the isolated point $[e]$ cannot belong to $\Lambda_X$.
If $\Aut(X)$ is infinite elementary, the stabilizer of every point of $\Lambda_X$ is infinite.
On the other hand, $\Aut(X,\pi_e)$ is finite, since $\MW(\pi_e)$ has finite index in $\Aut(X,\pi_e)$ and $\rk\MW(\pi_e)=r(e)=0$.
Since $\Aut(X,\pi_e)$ is the stabilizer of $[e]$, we have $[e]\notin\Lambda_X$.
If $\Aut(X)$ is finite, then the limit set is empty. 
Hence $E_X\cap\Lambda_X=\varnothing$. %in all cases.

\noindent(iii)
To treat the irrational case, %remaining case, 
we first show that $\Lambda_X=\Lambda_{A_X}$.
%We claim that $\Lambda_X=\Lambda_{A_X}$.
If $X$ satisfies \eqref{eq:char}, this follows from
Theorem~\ref{thm:cone-comparison}.
Otherwise, $X$ is supersingular in characteristic $2$.
By Lemma~\ref{lem:supersingular-mw} below,
there is a primitive nef isotropic class $e_0$ with $r(e_0)>0$.
For every $a\in A_X$, the class $ae_0$ is again primitive,
nef and isotropic, and $M_{ae_0}\cong M_{e_0}$.
Hence $r(ae_0)>0$ and
$
 A_X[e_0]\subseteq\Lambda_X\subseteq\Lambda_{A_X}.
$
If $A_X$ is non-elementary, its action on its limit set is
minimal, so
$\Lambda_{A_X}=\overline{A_X[e_0]}\subseteq\Lambda_X$,
as $\Lambda_X$ is closed.
If $A_X$ is elementary, the parabolic element fixing $[e_0]$
implies that $\Lambda_{A_X}=\{[e_0]\}$ and
$A_X[e_0]=\{[e_0]\}$.
Thus, the claim follows in either case.

%We next show that every irrational point of $\Lambda_X^{\mathrm{nef}}$ belongs to $\Lambda_{A_X}$, 
We follow the lattice-theoretic argument of \cite[Theorem~5.28]{Tak25}.
The group $G$ is a lattice when $\rho(X)\geq3$.
Choose $o\in\Amp(X)\cap\HH_X$ with trivial $G$-stabilizer,
and let $P$ be the closed Dirichlet domain for $G$ centered at $o$: 
\[
P:=\left\{x\in\HH_X\ \middle|\
d_{\HH_X}(x,o)\leq d_{\HH_X}(x,go)
\text{ for every }g\in G\right\}.
\]
The chamber decomposition gives
$P\subseteq D$ and $A_XP=D$; 
see \cite[Proposition~5.7 and Lemma~5.9]{Tak25}.
These assertions use only the lattice $L$ and its root reflections.

There is a $G$-invariant union $\mathcal U$ of pairwise disjoint
open horoballs based at rational isotropic rays such that
$(\HH_X\setminus\mathcal U)/G$ is compact;
see \cite[Theorem~3.31 and Lemma~5.8]{Tak25}.
We take $\mathcal U=\varnothing$ when $G$ is cocompact.
The set $P_0=P\setminus\mathcal U$ is compact.
Indeed, choose a compact set $Q\subseteq\HH_X\setminus\mathcal U$
with $GQ=\HH_X\setminus\mathcal U$.
If $x=gq\in P_0$, where $q\in Q$, the Dirichlet inequality gives
\[
 d(x,o)\leq d(x,go)=d(q,o).
\]
Thus $P_0$ is closed and bounded.
Moreover, $D\setminus\mathcal U=A_XP_0$.

Let $\xi\in\Lambda_X^{\mathrm{nef}}$ be irrational.
Since $D$ is convex, the geodesic ray from $o$ to $\xi$
lies in $D$.
It cannot eventually lie in $\mathcal U$:
otherwise, a tail of the ray would lie in a single horoball,
as the horoballs are pairwise disjoint.
Then $\xi$ would be the rational center of that horoball,
contradicting its irrationality.
Thus the ray meets $D\setminus\mathcal U$
arbitrarily far from $o$.
%The geodesic ray from $o$ to $\xi$ lies in $D$ and meets $D\setminus\mathcal U$ arbitrarily far from $o$.
%Otherwise, a tail of the ray would lie in one of the disjoint horoballs and $\xi$ would be its rational center.
We can therefore choose points $x_j\to\xi$ on the ray and
$a_j\in A_X$ such that $x_j\in a_jP_0$.
%rayが途中でhoroballに入ることは可能，だから$\mathcal Uに入らないような点列をchooseする
Compactness of $P_0$ gives a uniform bound on $d(x_j,a_jo)$.
Consequently $a_jo\to\xi$, proving that
$\xi\in\Lambda_{A_X}=\Lambda_X$.
Together with the rational cases (i)(ii), this completes the proof.
\end{pf}

\subsection{Quasi-isometries}\label{subsec:qi}

Let $G$ be a group with a finite symmetric generating set $S$.
The \emph{Cayley graph} $\Cay(G,S)$ has vertex set $G$ and an edge of
length one from $g$ to $gs$ for each $g\in G$ and $s\in S$.
We equip it with the path metric. Its restriction to the vertex set is
the \emph{word metric} on $G$
\[
 d_S(g,h)=|g^{-1}h|_S,
\]
where $|a|_S$ is the least length of a word in $S$ representing $a$.

\begin{definition}\label{def:qi}
Let $(Y,d_Y)$ and $(Z,d_Z)$ be metric spaces. A map $f:Y\to Z$ is a
\emph{quasi-isometric embedding} if there exist $\lambda\geq1$ and
$c\geq0$ such that
\begin{equation}\label{eq:qi}
 \lambda^{-1}d_Y(x,y)-c
 \leq d_Z(f(x),f(y))
 \leq\lambda d_Y(x,y)+c,
 \qquad x,y\in Y.
\end{equation}
It is a \emph{quasi-isometry} if, in addition, there exists $R\geq0$
such that, for every $z\in Z$, there exists $y\in Y$ with $d_Z(z,f(y))\leq R$. 
%such that every point of $Z$ is at distance at most $R$ from $f(Y)$.
We write $Y\qi Z$ if such a map exists.
\end{definition}

Quasi-isometry is an equivalence relation.
The identity between two word metrics on a finitely generated group is bi-Lipschitz 
(i.e. \eqref{eq:qi} holds with $c=0$).
Thus the quasi-isometry class of $G$ does not depend on $S$. Passing to a subgroup
of finite index or taking a quotient by a finite normal subgroup
preserves this class. In particular, commensurable finitely generated groups are quasi-isometric. 

The %following
\v Svarc--Milnor lemma is fundamental in relating
isometric actions to quasi-isometries:
if a group acts properly and cocompactly by isometries on a proper
geodesic metric space, then an orbit map is a quasi-isometry
\cite[Chapter~I.8, Proposition 8.19]{BH99}.

\begin{definition}\label{def:root-preserving}
Let $L$ and $L'$ be even hyperbolic lattices. An isometry of rational
quadratic spaces $\phi:L_\QQ\to L'_\QQ$ is \emph{root-preserving} if
\[
 \phi(\Delta(L))=\Delta(L'),%\qquad
 ~\text{ where }~
 \Delta(L):=\{r\in L\mid r^2=-2\}.
\]
%Here the roots are integral vectors in the specified lattices.
\end{definition}

\begin{proposition}\label{prop:root-qi}
Let $X$ and $Y$ be K3 surfaces, possibly over different algebraically
closed fields, satisfying \eqref{eq:char}. If there is a root-preserving rational isometry
\[
 \phi:\NS(X)_\QQ\xrightarrow{\sim}\NS(Y)_\QQ,
\]
then $\Aut(X)\qi\Aut(Y)$.
\end{proposition}

\begin{proof}
%Identify the two rational quadratic spaces by $\phi$, and 
Put $L_1:=\phi(\NS(X))$ and $L_2:=\NS(Y)$. 
These are commensurable in $\NS(Y)_\QQ$.
Their integral orthogonal groups $\Orth(L_1),\Orth(L_2)$ are commensurable in the rational
orthogonal group. To see this, choose an integer $a>0$ such that
$aL_1\subset L_2\subset a^{-1}L_1$.
The stabilizer of $L_2$ in $\Orth(L_1)$ has finite index, since
$\Orth(L_1)$ acts on the finite set of intermediate lattices between
$aL_1$ and $a^{-1}L_1$. Interchanging $L_1$ and $L_2$ proves the claim.

Replacing $\phi$ by $-\phi$ if necessary, we may assume that
$\phi(\cC_X)=\cC_Y$. Since $\phi$ identifies the root hyperplane
arrangements, $\phi(\Amp(X))$ is a chamber for $Y$.
Composing with an element of $W_Y$, we may further assume that
$\phi(\Amp(X))=\Amp(Y)$.
%この時点で\phiは(effはampleとの交点が正で特徴づけられるので)eff rootの集合，((-2)curve classはAmpのwallを定めるので)(-2)curve classの集合を保つことが分かる
%
Intersecting the group $\Orth(L_i)$ with the stabilizer of the common ample chamber shows that $\phi A_X\phi^{-1}$ and $A_Y$ are commensurable.
%commensurableを部分群に制限してもcommensurable
By Theorem~\ref{thm:cone-comparison}, their subgroups
$\phi\Gamma_X\phi^{-1}$ and $\Gamma_Y$ have finite index and are
therefore commensurable as well. The finite kernels of
$\Aut(X)\to\Gamma_X$ and $\Aut(Y)\to\Gamma_Y$ now give the assertion.
\end{proof}

\begin{comment}
The root-preserving condition concerns the integral root sets, and is
stronger than an isometry of rational quadratic spaces. It is precisely
what allows the ample chambers to be identified in this argument.
\end{comment}

\subsection{Relatively hyperbolic groups}\label{subsec:relhyp}

We recall some basic facts about hyperbolic and relatively
hyperbolic groups and their boundaries.
%We recall the boundary constructions used below. 
Standard references are \cite[Chapter~III.H]{BH99}, \cite{Bow12}, and \cite{GM08}.

A metric space is \emph{proper} if its closed bounded subsets are
compact, and \emph{geodesic} if any two points are joined by a path whose length equals their distance.
For subsets $B,C$ of a metric space $Y$, their Hausdorff distance is
\[
 d_{\mathrm H}(B,C)
 :=\inf\{R\geq0\mid B\subset N_R(C),\ C\subset N_R(B)\},
\]
where $N_R(B):=\{y\in Y\mid d_Y(y,B)\leq R\}$ is the $R$-neighbourhood of $B$, and the distance is $\infty$ if no such $R$ exists.

\begin{definition}\label{def:hyperbolic-space}
A geodesic metric space $Y$ is \emph{$\delta$-hyperbolic} if every
geodesic triangle is $\delta$-slim: each side lies in the
$\delta$-neighbourhood of the union of the other two sides.
It is \emph{hyperbolic} if it is $\delta$-hyperbolic for some
$\delta\geq0$.
\end{definition}

Hyperbolicity is invariant under quasi-isometry.
Let $Y$ be a proper geodesic hyperbolic space and fix $o\in Y$.
The \emph{Gromov boundary} $\partial_\infty Y$ is the set of geodesic
rays starting at $o$, modulo finite Hausdorff distance.
This definition is independent of $o$. 
For $x,y\in Y$, the Gromov product is
\[
 (x\mid y)_o
 :=\tfrac12\bigl(d_Y(o,x)+d_Y(o,y)-d_Y(x,y)\bigr)
 \in[0,\infty).
\]
For $\xi,\eta\in\partial_\infty Y$, we define
\[
 (\xi\mid\eta)_o
 :=\sup_{r,s}\lim_{t\to\infty}(r(t)\mid s(t))_o
 \in[0,\infty],
\]
where $r,s$ are geodesic rays from $o$ representing
$\xi,\eta$, respectively.
We equip $\partial_\infty Y$ with the visual topology,
characterised by $\xi_j\to\xi$ if and only if
$(\xi_j\mid\xi)_o\to\infty$.
With this topology, the boundary is compact.

\begin{definition}\label{def:visual}
A metric $d_{o,\epsilon}$ on $\partial_\infty Y$ is called
a \emph{visual metric} with parameter $\epsilon>0$
(based at $o$) 
if there exists $C\geq1$ such that
\begin{equation}\label{eq:visual}
 C^{-1}e^{-\epsilon(\xi\mid\eta)_o}
 \leq d_{o,\epsilon}(\xi,\eta)
 \leq Ce^{-\epsilon(\xi\mid\eta)_o}
 \qquad(\xi\ne\eta).
\end{equation}
\end{definition}

Visual metrics exist for all sufficiently small
$\epsilon>0$, and their induced topologies coincide
with the visual topology.

\begin{definition}\label{def:qs}
A homeomorphism $f:(Y,d_Y)\to(Z,d_Z)$ of metric spaces is \emph{quasisymmetric} 
if there exists an increasing homeomorphism $\eta:[0,\infty)\to[0,\infty)$ such
that
\[
 \frac{d_Z(f(x),f(a))}{d_Z(f(x),f(b))}
 \leq \eta\!\left(\frac{d_Y(x,a)}{d_Y(x,b)}\right)
\]
for all distinct $x,a,b\in Y$. We write $Y\qs Z$ if there is a
quasisymmetric homeomorphism between them.
\end{definition}

Every bi-Lipschitz homeomorphism is quasisymmetric.
Moreover, any two visual metrics on $\partial_\infty Y$
are quasisymmetrically equivalent.

A quasi-isometry between proper geodesic hyperbolic spaces
induces a homeomorphism of their Gromov boundaries,
which is quasisymmetric with respect to visual metrics
\cite[Proposition~6.3 and Theorem~6.5]{BS00}.
Moreover, quasi-isometries at bounded distance
induce the same boundary map
\cite[Proposition~6.3 (3)]{BS00}.

\begin{definition}\label{def:hyperbolic-group}
A finitely generated group $G$ is \emph{hyperbolic} 
if its %one, equivalently every, 
Cayley graph with respect to a finite generating set is a
hyperbolic space. For a finite generating set $S$, we write
$\partial G:=\partial_\infty\Cay(G,S)$,
called the \emph{Gromov boundary} of $G$.
\end{definition}

A hyperbolic group $G$ is \emph{elementary} if its Gromov boundary
$\partial G$ is finite (equivalently, if $G$ is finite or virtually infinite cyclic), and \emph{non-elementary} otherwise.
When $\Aut(X)$ is hyperbolic, this agrees with the notion of
elementarity introduced in Subsection~1.1, since
elementarity of $\Aut(X)$ is equivalent to virtual abelianity,
and hyperbolic groups contain no subgroup isomorphic to $\ZZ^2$.

Since Cayley graphs with respect to different finite generating sets are quasi-isometric and hyperbolicity of geodesic metric spaces is invariant under quasi-isometry, 
hyperbolicity of finitely generated groups is independent of the choice of a finite generating set and invariant under quasi-isometry.
The quasisymmetric class of the Gromov boundary $\partial G$
with a visual metric is also independent of the choice of a finite generating set.
%QIは境界のQSを誘導するので
%The Gromov boundary $\partial G$, equipped with a visual metric, is well-defined up to quasisymmetry, independently of the choice of a finite generating set.

\begin{thm}[Paulin \cite{Pau96}]\label{thm:qi-qs-hyperbolic}
Let $G$ and $H$ be hyperbolic groups.
Then $G$ and $H$ are quasi-isometric if and only if 
$\partial G$ and $\partial H$ are 
quasisymmetric with respect to visual metrics.
%quasisymmetrically homeomorphic with respect to visual metrics.
\end{thm}

%Thus the quasi-isometry classification of non-elementary hyperbolic groups is equivalent to the quasisymmetric classification of their Gromov boundaries.

We recall the definition of relative hyperbolicity in terms of cusped Cayley graphs.
For a connected graph $T$ with the path metric,
the \emph{combinatorial horoball} $\cH(T)$ associated with $T$
is the graph with vertex set $T^{(0)}\times\ZZ_{\geq0}$,
a copy of $T$ at level zero,
vertical edges from $(v,n)$ to $(v,n+1)$,
and horizontal edges from $(v,n)$ to $(w,n)$ whenever
$n\geq1$ and $0<d_T(v,w)\leq2^n$.
%Every edge has length one.

Let $G$ be finitely generated, and let $\cP$ be a finite collection of
infinite, finitely generated proper subgroups of $G$.
%proper subgpという仮定は意外と重要
Choose a finite generating set $S$ of $G$ such that
$S_P:=S\cap P$ generates $P$ for every $P\in\cP$.
For each left coset $gP$, attach a copy of
$\cH(\Cay(P,S_P))$ to $\Cay(G,S)$
by identifying its level-zero subgraph with the subgraph
of $\Cay(G,S)$ spanned by $gP$.
%by identifying $(p,0)$ with $gp$.
Equipping every edge with length one, the resulting graph
is a proper geodesic metric space.
We denote it by $\Cusp(G,\cP)$ and call it the
\emph{cusped Cayley graph} of the pair $(G,\cP)$.
Its quasi-isometry class is independent of the choice of $S$.
The cusped Cayley graph coincides with the cusped space
of Groves--Manning \cite[Definition~3.12]{GM08}.
%Definition~3.12は1-skeletonで定義してるからこの文章で問題ない

\begin{definition}[{Groves--Manning \cite{GM08}}]\label{def:relhyp}
The pair $(G,\cP)$ is \emph{relatively hyperbolic} if
its cusped Cayley graph $\Cusp(G,\cP)$ is hyperbolic. 
The subgroups in $\cP$ are called \emph{peripheral subgroups}. 
The Gromov boundary
$
 \partial(G,\cP):=\partial_\infty\Cusp(G,\cP)
$
is called the \emph{Bowditch boundary} of $(G,\cP)$.
%A visual metric on the Bowditch boundary is called a visual metric on $\partial(G,\cP)$.
\end{definition}

This definition is equivalent to several other definitions of
relative hyperbolicity \cite[Theorem~3.25]{GM08}.
We allow $\cP=\varnothing$; then $\Cusp(G,\varnothing)=\Cay(G,S)$, and relative hyperbolicity is ordinary hyperbolicity. 

Each combinatorial horoball attached to a left coset $gP$
determines a \emph{parabolic (fixed) point}
$\xi_{gP}\in\partial(G,\cP)$,
whose stabilizer in $G$ is $gPg^{-1}$. %i.e. conjugate of a peripheral subgroup

A relatively hyperbolic pair $(G,\cP)$ is \emph{elementary}
if its Bowditch boundary is finite (equivalently, if $G$ is
finite or virtually infinite cyclic), 
%定義にあるperipheral subgpがproperという条件がここで効いてる
and \emph{non-elementary} otherwise.
Equivalently, $(G,\cP)$ is non-elementary if and only if
$G$ contains a subgroup isomorphic to $F_2$.
For $G=\Aut(X)$, this coincides with elementarity of $\Aut(X)$, 
since 
%the latter is equivalent to virtual abelianity, and 
%elementarity of $\Aut(X)$ is equivalent to virtual abelianity, and
virtually abelian groups of rank at least two are not
relatively hyperbolic with respect to proper peripheral subgroups.

As in the case of hyperbolic groups, relative hyperbolicity
of a pair $(G,\cP)$ is independent of the choice of a finite
generating set.
The same holds for the quasisymmetric class of its Bowditch
boundary $\partial(G,\cP)$ equipped with a visual metric.
Moreover, relative hyperbolicity is invariant under
quasi-isometry \cite[Theorem~1.2]{Dru09}.
However, a quasi-isometry between two relatively hyperbolic
groups need not coarsely preserve their given peripheral families, and hence need not induce a quasisymmetric
homeomorphism between their Bowditch boundaries.
This motivates the following definition.

%An arbitrary hyperbolic model for the same topological boundary need not give this particular quasisymmetric class.

%For the metric assertion we use the cusped Cayley graphs above; it also holds for cusped spaces with constant horospherical distortion in the sense of \cite{HH21}.
%\cite{HH21}ではweak cusped spaceも導入している

\begin{comment}
Let $G$ be a group. 
For a family $\cP$ of subgroups of $G$, write
\[
G/\cP:=\{gP\mid g\in G,\ P\in\cP\}
\]
for the collection of left cosets.
\end{comment}

\begin{definition}\label{def:qi-pairs}
%\cite{HMPS23}と同じ定義
%\cite{HH21}よりも強い定義
Let $G,H$ be finitely generated groups, and let
$\cP,\cQ$ be finite families of subgroups of
$G,H$, respectively.
A quasi-isometry $f:G\to H$ is a
\emph{quasi-isometry of pairs}
$f:(G,\cP)\to(H,\cQ)$ if there is $D\geq0$ such that
\begin{enumerate}
\item for every $A\in G/\cP$, there is $B\in H/\cQ$ with
$d_{\mathrm H}(f(A),B)\leq D$;
\item for every $B\in H/\cQ$, there is $A\in G/\cP$ with
$d_{\mathrm H}(f(A),B)\leq D$.
\end{enumerate}
Such a quasi-isometry is also said to
\emph{coarsely preserve the peripheral families}
$\cP$ and $\cQ$.
\end{definition}

A family $\cP$ of subgroups of $G$ is \emph{almost malnormal} 
%\cite{HMPS23}と同じ定義
if
\[
 |P\cap gQg^{-1}|=\infty
 \quad\Longrightarrow\quad P=Q\ \text{and}\ g\in P
 \qquad(P,Q\in\cP,\ g\in G).
\]
Peripheral families of relatively hyperbolic groups are almost malnormal \cite[Proposition~2.36]{Osi06}. 
%, and their members are undistorted \cite[Lemma~5.4]{Osi06}.

\begin{theorem}\label{thm:relative-qi}
Let $G,H$ be finitely generated groups, and let $\cP,\cQ$ be finite
almost malnormal families of infinite finitely generated proper subgroups. 
%ここでalmost malnormalという仮定は必要
Suppose that $f:(G,\cP)\to(H,\cQ)$ is a quasi-isometry of pairs.
\begin{enumerate}
\item
The pair $(G,\cP)$ is relatively hyperbolic if and only if
$(H,\cQ)$ is relatively hyperbolic 
\textup{(}\cite[Theorem~C]{HMPS23} and \cite[Theorem~3.25]{GM08}\textup{)}. 
\item 
If these pairs are relatively hyperbolic, $f$ extends to a
quasi-isometry of their cusped Cayley graphs and induces a
quasisymmetric homeomorphism
\[
 \partial f:\partial(G,\cP)\xrightarrow{\sim}\partial(H,\cQ)
\]
for visual metrics.
This map sends parabolic points to parabolic points
\textup{(}\cite[Theorem~1.2]{HH21}, see also \cite[Corollary~1.2]{MS24}\textup{)}.
\end{enumerate}
\end{theorem}

%Part~\textup{(i)} follows from quasi-isometry invariance for pairs; see \cite[Theorem~C]{HMPS23} and \cite[Theorem~3.25]{GM08}. 
%Almost malnormality ensures that no refinement of the specified peripheral families is needed.
%Part~\textup{(ii)} is the cusp-extension theorem \cite[Theorem~1.2]{HH21}; see also \cite[Corollary~1.2]{MS24}.

%The boundary map is a map of metric spaces and is not asserted to be equivariant with respect to a group isomorphism.
%そもそもQIが群構造を保たないので当たり前

For a relative analogue of Paulin's theorem (\cref{thm:qi-qs-hyperbolic}), see \cite[Corollary~1.3]{MS24}.

The following is a consequence of \cite[Theorem~4.1]{BDM09};
see also \cite[Introduction]{HH21}.

\begin{theorem}[Behrstock--Dru\c tu--Mosher \cite{BDM09}]
\label{thm:nrh-peripherals}
Let $(G,\cP)$ and $(H,\cQ)$ be relatively hyperbolic pairs.
Suppose that no peripheral subgroup of either pair is relatively hyperbolic \textup{(}e.g. every such subgroup is virtually $\ZZ^r$
for some $r\geq2$\textup{)}.
%Suppose that every peripheral subgroup of either pair is not relatively hyperbolic (e.g. virtually $\ZZ^r$ for some $r\geq2$).
Then every quasi-isometry $G\to H$ coarsely preserves
these peripheral families.
%In particular, it induces a quasisymmetric homeomorphism of their Bowditch boundaries.
\end{theorem}

\begin{comment}
\cite{BDM09}から従うということで定着しているようだが，
(Hausdorff距離に関するboundは近傍への包含よりも強い主張なので)
実際に証明しようとすると\cite[Theorem~4.1]{BDM09}に加えてもう少し議論が必要(下記参照)
This follows from \cite[Theorem~4.1]{BDM09}.
Applying the theorem to a quasi-isometry and its quasi-inverse, together with the bounded intersection property of distinct peripheral cosets, shows that the peripheral families are coarsely preserved.

Corresponding peripheral subgroups have the same rank.
(peripheral subgpはundistortedよりperipheral間でもQSになるから)
\end{comment}

\begin{comment}
The rank restriction is necessary: virtually cyclic peripheral
subgroups are not automatically preserved by arbitrary
quasi-isometries. For example, the free group of rank two has Cantor boundary, whereas its relative boundary with respect to the cyclic subgroup generated by a punctured-torus peripheral element is a circle \cite[Example~1.1]{MS24}. 

One may omit hyperbolic peripheral subgroups from a peripheral family using transitivity of relative hyperbolicity \cite[Theorem~2.40]{Osi06}, but this changes the boundary in general.
\end{comment}

\vspace{2mm}

We now turn to the case of K3 surfaces.
%We now apply these facts to K3 surfaces. 

\begin{theorem}[{Kikuta \cite[Corollary~4.1]{Kik26}, Takatsu \cite[{Theorem~5.14}]{Tak25} over $\CC$}]\label{thm:k3-relative}
Let $X$ be a K3 surface satisfying \eqref{eq:char}. Then $\Aut(X)$ is
either virtually abelian or non-elementary relatively hyperbolic.
In the latter case, a peripheral family $\cP_X$ can be chosen
as a finite set of representatives of the conjugacy classes
of infinite groups $\Aut(X,\pi)$, where $\pi$ runs through
the genus-one fibrations.
Each peripheral subgroup is virtually abelian
of rank $\rk\MW(\pi)$.
\end{theorem}

\begin{comment}
For $\charac k\ne2$, or for a non-supersingular surface in
characteristic $2$ (i.e. \eqref{eq:char}??), 
this is \cite[Corollary~4.1 and the discussion
preceding it]{Kik26}; see also \cite[Theorem~5.14]{Tak25} over $\CC$.
\end{comment}

\begin{comment}
The same argument applies to supersingular surfaces in characteristic
$2$ with $\sigma(X)<10$, by \cite[Theorem~2.2]{Bra18} and the general criterion \cite[Theorem~3.1]{Kik26}. The description of the peripheral
subgroups uses \cite[Propositions~2.10 and~4.4]{KTT26}.
We do not use this theorem for the remaining case
$\charac k=2$, $\sigma(X)=10$.
\end{comment}

The virtually infinite cyclic subgroups in the peripheral family of Theorem~\ref{thm:k3-relative} can be omitted to obtain
another relatively hyperbolic structure on $\Aut(X)$ by \cite[Corollary~1.14]{DS05} (cf. \cite[Corollary~3.5]{FOY26}).
The resulting Bowditch boundary may differ from that of the full peripheral family, 
and we will specify the peripheral family whenever a Bowditch boundary is used.

%When applying Theorem~\ref{thm:nrh-peripherals} to $\Aut(X)$ in Sections 3 and 4,  we use only the peripheral groups of rank at least two. 

To state the hyperbolicity criterion without a characteristic
restriction \eqref{eq:char}, we deduce the following consequence of preceding results.

\begin{lemma}\label{lem:supersingular-mw}
Let $X$ be a supersingular K3 surface in characteristic $p>0$.
Then $X$ admits an elliptic fibration $\pi$ such that
\[
\rk\MW(\pi)\geq
\begin{cases}
2 & \text{if }p=2,\\
4 & \text{if }p>2.
\end{cases}
\]
In particular, $\Aut(X)$ contains $\ZZ^2$.
\end{lemma}

\begin{proof}
Let $\sigma=\sigma(X)$ be the Artin invariant of $X$, 
and $Y$ be a supersingular K3 surface of Artin
invariant 1 in characteristic $p$. 
By \cite[Lemma~2.1]{Bra18}, there is a finite-index
isometric embedding
\[
\NS(X)\hookrightarrow \NS(Y).
\]
We choose an elliptic fibration on $Y$ as follows.
If $p=2$, by \cite[Section~6, No.~8, and
Section~8.8]{ES15}, $Y$ admits an elliptic fibration whose associated root lattice is $E_6^{\oplus3}$. 
If $p>2$, we may take $Y=\Km(E\times E)$ for a supersingular elliptic curve $E$.
The projection $E\times E\to E$ induces an elliptic fibration on $Y$ with 4 singular fibers of type $I_0^*$,
whose associated root lattice is $D_4^{\oplus4}$; see \cite[Example~12.6 and Section~13.11]{SS10}.

Let $f\in\NS(Y)$ be the fiber class of the chosen fibration. 
%K3なのでfはprimitive
Then
\[
\rk\bigl(f^\perp_{\NS(Y)}/\ZZ f\bigr)_{\mathrm{rt}}
=
\begin{cases}
18 & \text{if }p=2,\\
16 & \text{if }p>2.
\end{cases}
\]

%Identify $\NS(X)$ with a finite-index sublattice of $\NS(Y)$, and 
Let $e$ be the primitive
generator of $\NS(X)\cap\ZZ_{\geq0}f$.
The inclusion induces a finite-index isometric
embedding
\[
M_e:=e^\perp_{\NS(X)}/\ZZ e
\lhook\joinrel\longrightarrow
f^\perp_{\NS(Y)}/\ZZ f:=M_f.
\]
Indeed, it is injective since $\ZZ e=\NS(X)\cap\ZZ f$ and 
$e^\perp_{\NS(X)}=\NS(X)\cap f^\perp_{\NS(Y)}$.
Every root of $M_e$ maps to a root of $M_f$. Consequently,
\[
\rk\bigl(M_e\bigr)_{\mathrm{rt}}
\leq
\begin{cases}
18 & \text{if }p=2,\\

16 & \text{if }p>2.
\end{cases}
\]

After changing sign and applying an element
of the Weyl group, we may assume that $e$
is nef. It defines a genus-one fibration
$\pi$ on $X$. By the Shioda--Tate formula \eqref{eq:mw-rank}, 
\[
\begin{aligned}
\rk\MW(\pi)
&=22-2-
\rk\bigl(M_e\bigr)_{\mathrm{rt}}\\
&\geq
\begin{cases}
2 & \text{if }p=2,\\
4 & \text{if }p>2.
\end{cases}
\end{aligned}
\]
Since a quasi-elliptic fibration has
Mordell--Weil rank zero, $\pi$ is elliptic.
Its Mordell--Weil group contains $\ZZ^2$,
as required.
\end{proof}

The following hyperbolicity criterion is a straightforward consequence of the preceding results and holds in arbitrary characteristic.
Over $\CC$, related criteria for non-elementary hyperbolicity were given in \cite[Theorems~1.1 and~1.3]{FOY26}. %by Fujiwara--Oguiso--Yu \cite[Theorems~1.1 and~1.3]{FOY26}.

%An alternative algebro-geometric proof that the existence of a subgroup isomorphic to $\ZZ^2$ forces a genus-one fibration of Mordell--Weil rank at least two can be found in \cite[Proposition~4.1]{FOY26}.

\begin{theorem}\label{thm:k3-hyperbolic}
Let $X$ be a K3 surface. %in arbitrary characteristic. 
The following are equivalent:
\begin{enumerate}
\item $\Aut(X)$ is hyperbolic;
\item $\Aut(X)$ contains no subgroup isomorphic to $\ZZ^2$;
\item either $X$ admits no genus-one fibration, or every
genus-one fibration $\pi$ on $X$ satisfies
$\rk\MW(\pi)\leq1$.
\end{enumerate}
%Moreover, $\Aut(X)$ is non-elementary hyperbolic if and only if it contains a nonabelian free group and contains no $\ZZ^2$.
\end{theorem}

\begin{proof}
A hyperbolic group contains no $\ZZ^2$, so \textup{(i)}
implies \textup{(ii)}. The translation action of $\MW(\pi)$
shows that \textup{(ii)} implies \textup{(iii)}.

Suppose that \textup{(iii)} holds.
By Lemma~\ref{lem:supersingular-mw}, $X$ satisfies
\eqref{eq:char}, so Theorem~\ref{thm:k3-relative} applies.
If $\Aut(X)$ is non-elementary, its peripheral subgroups
are virtually infinite cyclic by \textup{(iii)}.
Hence $\Aut(X)$ is hyperbolic by
\cite[Corollary~1.14]{DS05}, which proves \textup{(i)}.
It remains to consider the elementary case.
If $\Aut(X)$ is not virtually cyclic, 
%よって$\Aut(X)$はzero-entropy
then it preserves a genus-one fibration $\pi$ by \cite[Theorem~5.2]{BM23}
(cf.~\cite[Theorem~1.4]{Ogu07} over $\CC$).
Since $\MW(\pi)$ has finite index in $\Aut(X,\pi)=\Aut(X)$, 
we have $\rk\MW(\pi)\geq2$, contradicting \textup{(iii)}.
Thus $\Aut(X)$ is virtually cyclic and hence hyperbolic,
proving \textup{(i)}.
\end{proof}

\begin{comment}
Over $\CC$, the non-elementary criterion is
\cite[Theorem~1.1]{FOY26}. Fujiwara--Oguiso--Yu also give its
algebro-geometric form for $\rho(X)\geq5$: non-elementary
hyperbolicity is equivalent to the bound on Mordell--Weil ranks
together with the existence of an automorphism whose action on
$\NS(X)_\RR$ has spectral radius greater than one
\cite[Theorem~1.3]{FOY26}.

The proof above explains why the first criterion holds over any
algebraically closed field, including the characteristic-$2$ case
not covered by Theorem~\ref{thm:k3-relative}.
\end{comment}

\subsection{\texorpdfstring{\sier\ spaces}{Sierpinski spaces}}\label{subsec:sierpinski}

Let
$\SSph^m:=\{x\in\RR^{m+1}\mid\|x\|=1\}$,
where $\|\cdot\|$ is the standard Euclidean norm.
We equip $\SSph^m$ with its \emph{spherical metric}
\[
d_{\mathrm{sph}}(x,y):=\arccos\langle x,y\rangle,
\]
where $\langle\cdot,\cdot\rangle$ is the standard Euclidean inner product.
The \emph{chordal metric} on $\SSph^m$ is defined by
$d_{\mathrm{ch}}(x,y):=\|x-y\|$.
For the Poincar\'e ball model of $\HH^{m+1}$
with base point $o$ at its centre, we have
\begin{equation}\label{eq:spherical-visual}
e^{-(\xi\mid\eta)_o}
=\sin\!\left(\frac{d_{\mathrm{sph}}(\xi,\eta)}2\right)
=\frac12d_{\mathrm{ch}}(\xi,\eta).
\end{equation}
Thus the spherical metric, the chordal metric, and
any visual metric $d_{o,\epsilon}$ on $\partial_\infty \HH^{m+1}=\SSph^m$ with admissible
parameter $\epsilon>0$ are quasisymmetrically equivalent.
%\epsilon=1ならbi-Lipschitz/ 0<\epsilon\le1ならbi-LipschitzでないQS/ \epsilon=1ならそもそもmetricにならない

A \emph{round ball} is an open ball of radius less than
$\pi$ with respect to $d_{\mathrm{sph}}$;
its boundary is a \emph{round sphere}.

We also equip $\SSph^m$ with its \emph{spherical measure} $\mu_m$,
namely the Riemannian volume measure of the unit sphere,
or equivalently its $m$-dimensional Hausdorff measure
associated with $d_{\mathrm{sph}}$ with the standard normalization.
This measure is finite, and every nonempty open subset
of $\SSph^m$ has positive measure. %明らかにfull supportと同値
Moreover, bi-Lipschitz self-homeomorphisms of $\SSph^m$,
including the \mobius\ transformations introduced below,
preserve sets of spherical measure zero.

For a topological space $Z$, its \emph{topological dimension}
$\dim Z$ means its covering dimension: $\dim Z\leq n$ if every
finite open cover has an open refinement in which each point belongs
to at most $n+1$ members.
This is a topological invariant. %, and is distinct from Hausdorff dimension for a specified metric.

\begin{definition}
[{\cite{Can73}, see also \cite[Definition~1.1]{NW20}}]
\label{def:sierpinski}
%[{\cite[p.~107]{Can73}, see also \cite[Definition~1.1]{NW20}}]
For $m\geq0$, an \emph{$n$-dimensional \sier\ space}
is a subset $K\subset\SSph^{m+1}$ whose complementary components %$U_1,U_2,\ldots$ 
$\{U_i\}_{i}$ satisfy the following conditions. 
\begin{enumerate}
\item $\SSph^{m+1}\setminus U_i$ is homeomorphic to a
closed $(m+1)$-ball for every $i$;
\item $\overline U_i\cap\overline U_j=\varnothing$
for $i\ne j$;
%これはCannonの定義
%ambient homeoで移り合うという，埋込に関するより強い条件もある(このときperipheral domainをtame ballと呼ぶ)，この論文で扱うSierpinski spaceに関してはこの定義で問題ないことを全て確認済み．
%round ballはtame
\item $\bigcup_i U_i$ is dense in $\SSph^{m+1}$;
\item $\diam(U_i)\to0$ as $i\to\infty$ with respect to the spherical metric on $\SSph^{m+1}$.
\end{enumerate}
The sets $U_i$ and $\partial U_i$ are called its
\emph{peripheral domains} and \emph{peripheral spheres},
respectively. 
A \sier\ space is \emph{round} if every $U_i$ is a round ball.
\end{definition}

%We use the closed-complement condition in \textup{(i)}, as in \cite[Definition~1.1]{NW20}. 
%In higher dimensions, requiring the $U_i$ to be open topological balls (tame ball) instead would give a stronger condition, since their embeddings in our definition need not be tame.

The conditions in \cref{def:sierpinski} imply that $K$ is compact, that its peripheral domains form a countably infinite family, and
that $K$ is connected for $m\geq1$.
Round balls satisfy \textup{(i)} automatically.
An $m$-dimensional \sier\ space has topological dimension $m$.
For $m=0$ it is called a \emph{Cantor set}, and for $m=1$ it is called a \emph{\sier\ carpet}.

\begin{theorem}
[{Whyburn--Cannon \cite{Why58, Can73}, see also \cite[Theorem 1.2]{NW20}}]
\label{thm:sierpinski-unique}
For each $m\geq0$, any two $m$-dimensional \sier\ spaces are
homeomorphic.
\end{theorem}

\begin{comment}
For $m=0$, this is the characterisation of the Cantor set as a
nonempty compact metrizable perfect totally disconnected space.
For $m=1$, see \cite{Why58}; for higher dimensions, see \cite{Can73}.
Cannon's original statement excludes ambient dimension 4.
The 4-dimensional annulus theorem removes this restriction;
the resulting statement in all dimensions is recorded in
\cite[Theorem~1.2 and the paragraph preceding it]{NW20}.
\end{comment}

Although any two $m$-dimensional \sier\ spaces are
homeomorphic, they need not be quasisymmetrically equivalent
with respect to the spherical metrics.
This distinction is essential for the quasi-isometry classification in this paper.

\begin{definition}[{\cite[Section 1]{BKM09}}]\label{def:schottky}
A \emph{Schottky set} in $\SSph^m$ is a subset of the form
\[
S=\SSph^m\setminus\bigcup_{i\in I}B_i,
\]
where the $B_i$ are pairwise disjoint round open balls and
$|I|\geq3$.
It is equipped with the metric induced by the spherical metric.
\end{definition}

%This is the convention of \cite{BKM09}. 

Every round \sier\ space is a Schottky set.
However, the disjointness of the open balls in
\cref{def:schottky} allows their boundaries to intersect.
Thus a Schottky set need not be a \sier\ space.

\begin{comment}
It is a $(m-1)$-dimensional round \sier\ space if the closures of its complementary balls are pairwise disjoint and their union is dense. 
For infinitely many disjoint round balls, the null condition $\diam(B_i)\to0$ is automatic, since the sphere has finite measure and balls with a fixed positive diameter have measure bounded below.
\end{comment}

The group of \mobius\ transformations of $\SSph^m$, including
orientation-reversing transformations, is identified with
$\Isom(\HH^{m+1})$ through the action on the ideal boundary. 
%Equivalently, it is generated by reflections in round spheres.
Two subsets $S,S'\subset\SSph^m$ are \emph{\mobius\ equivalent} if
$M(S)=S'$ for some \mobius\ transformation $M:\SSph^m\to\SSph^m$.
Every \mobius\ transformation is quasisymmetric for the spherical
metric and preserves the class of round balls.
%Mobius変換はbi-Lipschitz homeoなのでQIかつQS/ しかし一般にisometryではない/ open ballのroundnessは保つが半径は変わる/ spherical measure zero setを保つ

\begin{theorem}[{Bonk--Kleiner--Merenkov~\cite[Theorem~1.1]{BKM09}}]\label{thm:bkm}
Let $m\geq2$, and let $S,S'\subset\SSph^m$ be Schottky sets.
If $\mu_m(S)=0$, then every quasisymmetric homeomorphism
$f:S\to S'$ is the restriction of a \mobius\ transformation of
$\SSph^m$.
\end{theorem}
Consequently, for a Schottky set
$S\subset\SSph^m$ of spherical measure zero and another Schottky set
$S'\subset\SSph^m$ with $m\geq2$,
we have
\begin{equation*}%\label{eq:qs-mobius}
 S\qs S'
 \quad\Longleftrightarrow\quad
 S\text{ and }S'\text{ are \mobius\ equivalent}.
\end{equation*}

\begin{comment}
The dimension and measure assumptions are essential to the stated
result. It does not apply to Cantor subsets of $\SSph^1$.
Moreover, topological dimension $m-1$ does not imply vanishing
$m$-dimensional spherical measure.
\end{comment}

%If a boundary in the preceding subsection is identified quasisymmetrically with a measure-zero Schottky set, then Theorem~\ref{thm:bkm} can be applied to its quasisymmetric boundary maps. 
%Such an identification is a metric assertion; a homeomorphism alone does not suffice.

We finish by recalling the role of \sier\ carpets among one dimensional boundaries of hyperbolic groups. 
A \emph{nontrivial splitting} of a group over a subgroup
is an amalgamated-product or HNN decomposition over that
subgroup whose Bass--Serre tree has no global fixed point.
A compact space is \emph{planar} if it embeds in $\SSph^2$.

\begin{theorem}[{Kapovich--Kleiner \cite[Theorem~4]{KK00}}]\label{thm:kk}
Let $G$ be a non-elementary hyperbolic group with
$\dim\partial G=1$. Suppose that $G$ admits no nontrivial splitting
over a finite or virtually cyclic subgroup. Then $\partial G$ is
homeomorphic to a circle, a \sier\ carpet, or the Menger curve.
If $\partial G$ is planar and is not a circle, then it is a
\sier\ carpet.
\end{theorem}

The last assertion of \cref{thm:kk} follows because the Menger curve is non-planar. 
The circle case, characterized in
\cref{thm:boundary-free-surface} (ii),
is not excluded by the no-splitting assumption.
Thus, among the remaining one-dimensional boundaries
under this assumption, the \sier\ carpet is precisely
the planar case.
In this sense, carpet boundaries form an irreducible case
in the quasi-isometry classification of hyperbolic groups
with one-dimensional boundaries.

\begin{comment}
The trichotomy is \cite[Theorem~4]{KK00}; 
In \cref{thm:kk}, the last assertion follows
because the Menger curve is non-planar. 
The circle alternative cannot be omitted solely from the no-splitting assumption.
\textcolor{red}{
See \cref{thm:boundary-free-surface} (ii) for the characterization of the circle case. 
}
Thus, after excluding the circle case, the \sier\ carpet is the planar one-dimensional boundary occurring for groups that cannot be further decomposed over finite or virtually cyclic subgroups.
This is the sense in which carpet boundaries form an irreducible case in the quasi-isometry classification.
\end{comment}

\begin{comment}
Theorem~\ref{thm:sierpinski-unique} fixes their topology, whereas their quasisymmetric structures retain the information relevant to
quasi-isometry. 
For measure-zero round realizations, this additional information is governed by \eqref{eq:qs-mobius}.
\end{comment}
% !TEX root = ../main.tex
% !TEX root = ../main.tex
% Integrated from section3_1_sierpinski.tex; full change log: ../CHANGELOG-ja.md.
\section{\texorpdfstring{K3 surfaces of \sier\ type}{K3 surfaces of Sierpinski type}}
\label{sec:sierpinski-type}

Throughout this section, let $X$ be a K3 surface over an algebraically closed field $k$ of arbitrary characteristic. 
We use the notation of Section~\ref{sec:preliminaries}. 
%In particular, $\NS(X)$, $\rho(X)=\rk \NS(X)$, and $\Gamma_X$ is the image of $\Aut(X)$ in $\Orth^+(\NS(X))$. 
%A \emph{$(-2)$-curve} is a smooth rational curve on $X$.
We denote the set of $(-2)$-curves (or equivalently, smooth rational curves) by $\cR_X$, and identify each such
curve with its class in $\NS(X)$ when taking intersection numbers.

\subsection{Definition and fundamental properties}\label{subsec:fundamental}

For each $C\in\cR_X$, define the open ball
\[
 B_C:=\{[v]\in\partial_\infty\HH_X\mid(v,C)<0\}.
\]
The description \eqref{eq:ample} of the ample cone in terms of $(-2)$-curves gives
\begin{equation}\label{eq:amp-boundary}
 \Lambda_X^{\mathrm{nef}}
 =\partial_\infty\HH_X\setminus\bigcup_{C\in\cR_X}B_C.
\end{equation}

\begin{lemma}\label{lem:root-ball-comparison}
Let $X$ be a K3 surface with $\rho(X)\geq 4$,
and let $C,D\in\cR_X$ be distinct.
Neither $B_C$ nor $B_D$ contains the other.
Moreover,
\begin{equation}\label{eq:root-ball-comparison}
\overline{B_C}\cap\overline{B_D}=\varnothing
\quad\Longleftrightarrow\quad
(C,D)\geq 3.
\end{equation}
\end{lemma}

\begin{proof}
Put $a=(C,D)\in\ZZ_{\geq 0}$.
If $a\geq 3$, then $(C+D)^2=2a-4>0$.
Since $C+D$ is effective, it belongs to
the positive cone $\cC_X$.
Thus $(v,C+D)>0$ for every nonzero
$v\in\overline{\cC_X}$, which excludes
$(v,C)\leq 0$ and $(v,D)\leq 0$
simultaneously. Hence
$\overline{B_C}\cap\overline{B_D}=\varnothing$.
If $a=2$, the vector $e=C+D$ is nonzero,
isotropic, and belongs to $\overline{\cC_X}$.
Since $(e,C)=(e,D)=0$, we obtain
$[e]\in\partial B_C\cap\partial B_D$.
Any isotropic vector orthogonal to both
$C$ and $D$ is proportional to $e$,
since $e^\perp/\RR e$ is negative definite.
Thus $[e]$ is the unique common boundary point.
Moreover, $(v,e)\geq 0$ for every
$v\in\overline{\cC_X}$, so
$(v,C)$ and $(v,D)$ cannot both be negative.
Hence $B_C$ and $B_D$ are externally tangent
at $[e]$.
If $a=0$ or $1$, the span of $C,D$ is
negative definite, so its orthogonal
complement has signature $(1,\rho-3)$.
As $\rho\geq 4$, this complement contains
a nonzero real isotropic vector in
$\overline{\cC_X}$, giving
$\partial B_C\cap\partial B_D\neq\varnothing$.
At any common boundary point $[v]$,
the two boundary spheres intersect
transversely: otherwise, some nonzero
linear combination of $C$ and $D$ would
be proportional to $v$, contradicting
the negative definiteness of their span.

Thus neither ball contains the other in any case:
their closures are disjoint when $a\geq 3$,
their boundaries are externally tangent
when $a=2$, and their boundaries intersect
transversely when $a=0$ or $1$.
This proves both assertions.
\end{proof}

\begin{definition}\label{def:sierpinski-type}
A K3 surface $X$ is \emph{of \sier\ type} 
if the following conditions hold:
\begin{enumerate}
\item $\rho(X)\geq4$; 
\item $X$ contains a $(-2)$-curve;
\item $(C,D)\geq3$ for any distinct $C,D\in\cR_X$.
\end{enumerate}
\end{definition}

\begin{remark}
We explain the restriction $\rho(X)\geq4$ in
\cref{def:sierpinski-type}.
For $\rho(X)\leq2$, the group $\Aut(X)$ is virtually
cyclic, and the ideal boundary $\partial_\infty\HH_X$
has at most two points.

When $\rho(X)=3$, the ideal boundary
$\partial_\infty\HH_X$ is a circle, and the ideal
boundary of a root hyperplane consists of two points.
In particular, the latter is 0-dimensional, as
are the isolated points of the nef boundary arising
from genus-one fibrations of Mordell--Weil rank 0.
Moreover, two distinct root hyperplanes may intersect
inside $\HH_X$ while their ideal boundaries are
disjoint. Indeed, for roots $r\ne\pm s$ with
$|(r,s)|<2$, their common orthogonal complement
has signature $(1,\rho(X)-3)$.
It contains no nonzero isotropic vector when
$\rho(X)=3$, whereas for $\rho(X)\geq4$
its isotropic rays form an
$\SSph^{\rho(X)-4}$ in the intersection of the
ideal boundaries.
Thus the equivalence in
\cref{lem:root-ball-comparison} fails
in Picard rank 3.

There are also obstacles to quasi-isometry rigidity (\cref{thm:rigidity-sierpinski}).
A 0-dimensional \sier\ space is a Cantor set
and hence totally disconnected, unlike the
positive-dimensional \sier\ spaces considered here.
Moreover, the rigidity theorem of
Bonk--Kleiner--Merenkov
\cite[Theorem~1.1]{BKM09} requires the ambient
sphere to have dimension at least 2.
In fact, non-elementary virtually free automorphism
groups occur in different Picard ranks, as illustrated
by \cref{ex:rank-three-free,ex:rank-four-free},
so their quasi-isometry classes need not determine
the rational N\'eron--Severi lattices.

Finally, in Picard rank 3, the intersection
condition in \cref{def:sierpinski-type} forces
every genus-one fibration to have Mordell--Weil
rank 1, by the Shioda--Tate formula.
Thus its infinite peripheral subgroups are
virtually cyclic, and arbitrary quasi-isometries
need not coarsely preserve the peripheral family
\cite[Example~1.1]{MS24}, see also \cref{thm:nrh-peripherals}.
These distinctions motivate the restriction
$\rho(X)\geq4$.
\end{remark}

The followings are fundamental properties of K3 surfaces of \sier\ type.

\begin{theorem}%[Theorem C]
\label{thm:fundamental}
Let $X$ be a K3 surface of \sier\ type. Then the following hold.
\begin{enumerate}
\item\label{item:rank} $\rho(X)\in\{4,5\}$.
\item\label{item:fibrations} Every genus-one fibration
$\pi:X\to\PP^1$ is elliptic and has no reducible fibers, i.e. its singular fibers are of type $I_1$ or $I\!I$, and
\[
 \MW(\pi)\cong\ZZ^{\rho(X)-2}.
\]
%Here $\MW(\pi)$ denotes the Mordell--Weil group of the Jacobian fibration, as in Section~\ref{sec:preliminaries}.
\item\label{item:limit} $\Lambda_X=\Lambda_X^{\mathrm{nef}}$.
\item\label{item:non-elementary} $\Aut(X)$ is non-elementary.
\item\label{item:curves} $X$ contains infinitely many $(-2)$-curves.
\item\label{item:hyperbolic} The following conditions are equivalent:
$\Aut(X)$ is hyperbolic; $\Aut(X)$ is convex-cocompact;
$X$ admits no genus-one fibration; and $\NS(X)$ is anisotropic.
If these conditions hold, then $\rho(X)=4$.
\end{enumerate}
\end{theorem}

\begin{comment}
\begin{remark}\label{rem:rank-three}
Let $X$ have Picard rank 3 and satisfy the two conditions in
Definition~\ref{def:sierpinski-type}. The same proof gives
\textup{\ref{item:fibrations}--\ref{item:curves}} of
Theorem~\ref{thm:fundamental}, with $\MW(\pi)\cong\ZZ$ whenever
a genus-one fibration exists. Moreover, $\Aut(X)$ is always
non-elementary hyperbolic, while $\Gamma_X$ is convex-cocompact if
and only if there is no genus-one fibration. Thus, if the rank
restriction is omitted from the definition, the corresponding forms
of \textup{\ref{item:rank}} and the hyperbolicity assertion in
\textup{\ref{item:hyperbolic}} are
\[
 3\leq\rho(X)\leq5,\qquad
 \Aut(X)\text{ is hyperbolic}
 \ \Longleftrightarrow\
 \rho(X)=3\ \text{or $X$ has no genus-one fibration}.
\]
\end{remark}
\end{comment}

\begin{proof}
We first show that every genus-one fibration has
irreducible fibers. Let $e$ be its fiber class.
If $C,D$ are distinct components of a reducible
fiber, then they are $(-2)$-curves and, since
$e^\perp$ is negative semidefinite,
\[
(C+D)^2=-4+2(C,D)\leq0.
\]
Thus $(C,D)\leq2$, contrary to the definition.

For \textup{\ref{item:rank}}, suppose that
$\rho(X)\geq6$, and choose $C\in\cR_X$.
The lattice $C^\perp_{\NS(X)}$ has signature
$(1,\rho(X)-2)$ and rank at least five.
By Meyer's theorem, it contains a nonzero
integral isotropic vector $e$.
After changing its sign and applying an element
of $W_X$, we may assume that $e$ is nef.
It defines a genus-one fibration.
By Riemann--Roch, one of $C$ and $-C$,
after applying the same Weyl group element,
is represented by an effective divisor
supported in fibers.
Since all fibers are irreducible, such a divisor
has square zero, a contradiction.
Hence $\rho(X)\leq5$.

For \textup{\ref{item:fibrations}}, the absence of
reducible fibers and the Shioda--Tate formula \eqref{eq:mw-rank} give
\[
\rk\MW(\pi)=\rho(X)-2>0.
\]
Since a quasi-elliptic fibration has Mordell--Weil
rank zero, $\pi$ is elliptic, and its singular
fibers are of type $I_1$ or $I\!I$.
The Jacobian fibration has the same fiber types
\cite[Chapter~11, \S4.1]{Huy16}.
In particular, it has no reducible fibers.
For any section $P\ne O$, the height formula gives
\[
\langle P,P\rangle=4+2(P,O)>0
\]
\cite[\S11.8]{SS10}.
Thus $\MW(\pi)$ is torsion-free, proving
\textup{\ref{item:fibrations}}.

For \textup{\ref{item:limit}}, \cref{limit-set-nef-boundary} shows that
$\Lambda_X^{\mathrm{nef}}\setminus\Lambda_X$
consists precisely of the nef isotropic rays corresponding
to genus-one fibrations with Mordell--Weil rank zero.
By \textup{\ref{item:fibrations}}, every genus-one fibration
has Mordell--Weil rank $\rho(X)-2>0$.
Hence $\Lambda_X=\Lambda_X^{\mathrm{nef}}$.

For \textup{\ref{item:non-elementary}}, choose
$C\in\cR_X$.
By Lemma~\ref{lem:root-ball-comparison},
the sphere $\partial B_C$ is contained in
$\Lambda_X^{\mathrm{nef}}=\Lambda_X$.
Since $\rho(X)\geq4$, this sphere is infinite.
Hence $\Aut(X)$ is non-elementary.

For \textup{\ref{item:curves}}, 
since $X$ contains a $(-2)$-curve and $\Aut(X)$ is infinite by \textup{\ref{item:non-elementary}}, 
the assertion follows from \cite[Chapter~8, Corollary~4.7]{Huy16}.
The same argument applies in positive characteristic by the cone conjecture 
(\cite[Theorem~6.1]{LM18} and
\cite[Theorem~3.1.4]{Takamatsu25}),
since $X$ is not supersingular by \textup{\ref{item:rank}}.

\begin{comment}
For \textup{\ref{item:curves}}, suppose that
$\cR_X$ is finite.
A finite-index subgroup $\Gamma_0$ of $\Gamma_X$
fixes the class of some $C\in\cR_X$.
Consequently,
\[
\Lambda_X=\Lambda_{\Gamma_0}
\subseteq
\partial_\infty(C^\perp\cap\HH_X)
\cong
\SSph^{\rho(X)-3}
\subset
\SSph^{\rho(X)-2},
\]
where $\Lambda_{\Gamma_0}$ is the limit set of $\Gamma_0$. 
On the other hand,
Lemma~\ref{lem:root-ball-comparison} and
\textup{\ref{item:limit}} imply that
\[
\Lambda_X
=\partial_\infty\HH_X
\setminus\bigcup_{D\in\cR_X}B_D
\]
contains a nonempty open subset of
$\partial_\infty\HH_X$.
This is a contradiction, since
$\partial_\infty(C^\perp\cap\HH_X)$
has empty interior.
Thus $\cR_X$ is infinite.
\end{comment}

Finally, by \textup{\ref{item:fibrations}} and
\cite[Theorem~4.11]{Kik26},
$\Aut(X)$ is convex-cocompact if and only if
$X$ admits no genus-one fibration.
In this case, $\Aut(X)$ is hyperbolic.
Conversely, if $X$ admits a genus-one fibration,
then $\Aut(X)$ contains a subgroup isomorphic to
$\ZZ^{\rho(X)-2}$, with $\rho(X)-2\geq2$,
and hence is not hyperbolic.
Moreover, 
%a nonzero primitive isotropic vector in $\NS(X)$ can be made nef by changing its sign and applying an element of $W_X$.
it follows that $X$ admits no genus-one fibration
if and only if $\NS(X)$ is anisotropic.
Meyer's theorem excludes anisotropy when
$\rho(X)=5$.
This proves \textup{\ref{item:hyperbolic}}.
\end{proof}

The definition can be expressed entirely in terms of the N\'eron--Severi lattice. 
Recall that $\Delta_X:=\{r\in \NS(X)\mid r^2=-2\}$.

\begin{lemma}\label{lem:numerical}
A K3 surface $X$ is of \sier\ type if and only if
the following conditions hold:
\begin{equation}\label{eq:numerical}
 \rho(X)\geq4,\qquad \Delta_X\ne\varnothing,\qquad
 |(r,s)|\geq3\quad
 \text{for all }r,s\in\Delta_X\text{ with }r\ne\pm s.
\end{equation}
%Under these conditions, $\Delta_X$ is infinite. (-2)曲線が無限個あるんだから当たり前
%Thus one may equivalently replace $\Delta_X\ne\varnothing$ by the requirement that $\NS(X)$ contain infinitely many roots.
\end{lemma}

\begin{proof}
Recall that the Riemann--Roch theorem implies that
every integral class of square at least $-2$ is effective up to sign. In particular,
$\Delta_X\ne\varnothing$ if and only if
$\cR_X\ne\varnothing$.

Suppose first that $X$ is of \sier\ type, and let
$r,s\in\Delta_X$ with $r\ne\pm s$.
If $|(r,s)|<2$, their span is negative definite.
Choose an integral positive vector $h$ orthogonal
to $r$ and $s$, and apply an element of $W_X$
to make $h$ nef.
Up to sign, the transformed roots are represented
by effective divisors supported on $(-2)$-curves
orthogonal to $h$.
There is at most one such curve, since
\[
(C+D)^2=-4+2(C,D)>0
\]
for distinct $(-2)$-curves $C,D$, contrary to
the negative definiteness of $h^\perp$.
Thus $r=\pm s$, a contradiction.
If $|(r,s)|=2$, then $e=r\pm s$ is a nonzero
isotropic vector orthogonal to $r$ for a suitable
choice of sign.
After making $e$ nef, the transformed root is effective up to sign and supported in fibers of the associated genus-one fibration.
This contradicts Theorem~\ref{thm:fundamental}~\textup{\ref{item:fibrations}}.
Hence $|(r,s)|\geq3$.

Conversely, suppose that \eqref{eq:numerical} holds.
Then $\cR_X\ne\varnothing$, and for distinct
$(-2)$-curves $C,D$, we have $(C,D)\geq3$. 
%we have $(C,D)\geq0$, so \eqref{eq:numerical} gives $(C,D)\geq3$.
Thus $X$ is of \sier\ type.
\end{proof}

A root-preserving rational isometry preserves \sier\ type.

\begin{cor}\label{prop:rigidity-general-direction}
Let $X_1$ be a K3 surface of \sier\ type and 
$X_2$ a K3 surface, possibly defined over a different
algebraically closed field.
If there is a root-preserving rational isometry
$\NS(X_1)_\QQ\xrightarrow{\sim}\NS(X_2)_\QQ$, then $X_2$ is of \sier\ type. %and $\Aut(X_1)\qi\Aut(X_2)$.
\end{cor}

\begin{proof}
A root-preserving rational isometry preserves the rank, the nonemptiness of the root set, and the inequalities
$|(r,s)|\geq3$ for roots $r\ne\pm s$.
The numerical characterization in Lemma~\ref{lem:numerical} therefore implies that $X_2$ is of \sier\ type.
%\cref{prop:root-qi} gives the quasi-isometry.
\end{proof}

We now give examples. 
All lattices below are even hyperbolic of rank 4 or 5, 
%All lattices below are even, have signature $(1,r-1)$ with $r=4$ or $5$, 
hence occur as the N\'eron--Severi lattices of complex K3 surfaces by the primitive embedding theorem and the surjectivity of the period map \cite[Corollary~14.3.1 and Theorem~6.3.1]{Huy16}. 
It is easy to verify the conditions of Lemma~\ref{lem:numerical} directly.
%We verify the conditions of Lemma~\ref{lem:numerical} directly.

\begin{example}\label{ex:Sierpinski-lattice}
The following lattices occur as the Néron–Severi lattices of complex K3 surfaces of Sierpiński type. 
\begin{enumerate}
\item
Picard rank 4, anisotropic case:
%Case of Picard rank 4, without isotropic vectors: 
\begin{equation*}%\label{eq:anisotropic-matrix}
L_1:=
 \begin{pmatrix}
 2&3&3&3\\
 3&-2&4&4\\
 3&4&-2&4\\
 3&4&4&-2
 \end{pmatrix}.
\end{equation*}

\item
Picard rank 4, isotropic case:
%Case of Picard rank 4, with isotropic vectors: 
\begin{equation*}%\label{eq:isotropic-four}
 L_2:=U(2)\oplus\langle-6\rangle^{\oplus2},
 \qquad
 U(2)=\begin{pmatrix}0&2\\2&0\end{pmatrix}.
\end{equation*}

\item
Picard rank 5:
%Case of Picard rank 5: 
\begin{equation*}%\label{eq:rank-five}
 L_3:=\langle4\rangle\oplus\langle-8\rangle
       \oplus\langle-12\rangle^{\oplus2}\oplus\langle-2\rangle.
\end{equation*}
\end{enumerate}
\end{example}

\begin{example}[Rational isometries preserving roots]
\label{ex:root-qi}
The N\'eron--Severi lattice need not be determined
by the quasi-isometry class of the automorphism group,
even among K3 surfaces of \sier\ type
(cf. \cref{thm:rigidity-sierpinski}).
Let $e,f,g_1,g_2$ be a basis of the lattice $L_2$
in \cref{ex:Sierpinski-lattice}~(ii), with
\[
(e,f)=2,\qquad g_1^2=g_2^2=-6,
\]
and all other pairings zero.
We identify $ae+bf+cg_1+dg_2$ with $(a,b,c,d)$
and define
\[
 L_4:=\{(a,b,c,d)\in L_2\mid a\equiv b\pmod 3\}. 
\]
Every root $(a,b,c,d)\in L_2$ satisfies
$
4ab-6(c^2+d^2)=-2,
$
hence $ab\equiv1\pmod3$ and thus $a\equiv b\pmod3$.
Consequently,
\[
 \Delta(L_4)=\Delta(L_2).
\]
The inclusion $L_4\subset L_2$ has index 3
and induces a root-preserving rational isometry.
\begin{comment}
In the basis $e+f,3f,g_1,g_2$, a Gram matrix of
$L_4$ is
\[
 \begin{pmatrix}
 4&6&0&0\\
 6&0&0&0\\
 0&0&-6&0\\
 0&0&0&-6
 \end{pmatrix}.
\]
\end{comment}
Since we have $\det L_4=9\det L_2$, 
\begin{comment}
On the other hand,
\[
 \det L^{(1)}=-144,\qquad \det L^{(3)}=-1296,
\]
so 
\end{comment}
these lattices are not integrally isometric.
Choose complex K3 surfaces $X_n$ with
$\NS(X_n)\cong L_n$ for $n=2,4$.
Both are of \sier\ type by \cref{lem:numerical},
and Proposition~\ref{prop:root-qi} gives
\[
 \Aut(X_2)\qi\Aut(X_4).
\]
%The notation $L^{(n)}$ here denotes a finite-index sublattice with the restricted form. It is different from the usual rescaling $L(n)$.
%For an even lattice $L$, the rescaled lattice $L(n)$ has no roots when $n>1$, and hence cannot itself be the N\'eron--Severi lattice of a K3 surface of \sier\ type.
\end{example}

%The preceding constructions use the period map and do not give defining equations for the surfaces. 
%This leads to the following problem.

The following problem is natural to consider.

\begin{problem}\label{prob:projective-models}
Construct explicit geometric models of K3 surfaces of \sier\ type.
%In particular, construct such models with N\'eron--Severi lattices $L_1$, $L_2$ and $L_3$ in Examples~\ref{ex:anisotropic}--\ref{ex:five}, and verify their N\'eron--Severi lattices.
\end{problem}
% !TEX root = ../main.tex
% Integrated from section3_2_round_limit_sets.tex; full change log: ../CHANGELOG-ja.md.
\subsection{Limit sets}
%\subsection{Round limit sets and a recognition question}
%\label{subsec:characterization}
\label{subsec:round-limit-sets}

We prove that a K3 surface of \sier\ type has a round \sier\ limit
set, whose complementary components correspond to its $(-2)$-curves.
%We also formulate the converse as a question and give a finiteness criterion over $\CC$.

%We use the notation of Section~\ref{sec:preliminaries} and Subsection~\ref{subsec:fundamental}. 
%Unless otherwise stated, $X$ is a K3 surface over an algebraically closed field of arbitrary characteristic, $\rho(X)\geq4$, and $m=\rho-2$.

An isometry from $\HH_X$ to the Poincar\'e ball model of $\HH^{\rho(X)-1}$ induces a quasisymmetric homeomorphism
\[
\partial_\infty\HH_X\qs\SSph^{\rho(X)-2}
\]
with respect to a visual metric on $\partial_\infty\HH_X$ and the spherical metric on $\SSph^{\rho(X)-2}$ by \eqref{eq:spherical-visual}. 
Under this identification, each $B_C~(C\in\cR_X)$ is a round open ball, since the defining inequality $(v,C)<0$ describes a spherical cap in Lorentz coordinates.
Any two such identifications differ by a \mobius\ transformation of $\SSph^{\rho(X)-2}$, so roundness of $B_C$ is independent of the chosen isometry.
See also \cite[Chapter I.6, \S6.7]{BH99}. 

%We use the spherical metric and the spherical measure on $\SSph^m$.

\begin{lem}\label{proper-lim-measure-zero}
Let $X$ be a K3 surface satisfying \eqref{eq:char}.
If $\cR_X\ne\varnothing$, then
$\Lambda_X\subset\SSph^{\rho(X)-2}$
has zero spherical measure.
\end{lem}

\begin{pf}
The action of $\Aut(X)$ on $\HH_X$ is geometrically finite by \cite[Theorem~A]{Kik26}. 
Since we have $\Lambda_X\subseteq\Lambda_X^{\mathrm{nef}}$,
\eqref{eq:amp-boundary} and $\cR_X\ne\varnothing$
imply that $\Lambda_X$ is a proper closed subset of $\SSph^{\rho(X)-2}$, 
%proper = 真部分集合の意味
%closedなので(補集合がopenなので)full spherical measureを持たない
hence $\Lambda_X$ has zero spherical measure by Ahlfors' theorem \cite[Theorem~8.4.2]{Thu22}. 
%for geometrically finite Kleinian groups 
\end{pf}

\begin{theorem}%[Theorem B]
\label{thm:round-limit-set}
If $X$ is of \sier\ type, then $\Lambda_X$ is a
$(\rho(X)-3)$-dimensional round \sier\ space.
Each complementary component is precisely the
round open ball $B_C~(C\in\cR_X)$, 
and the assignment $C\mapsto B_C$ is a bijection.
%\[
% B_C=\{[v]\in\partial_\infty\HH_X\mid(v,C)<0\},
% \qquad C\in\cR_X.
%\]
\end{theorem}

\begin{proof}
We verify the conditions in Definition~\ref{def:sierpinski}.
%
%Set $m:=\rho(X)-2$.
By \cref{thm:fundamental} (i)(iii)(v),
we have $\rho(X)\in\{4,5\}$, $\Lambda_X=\Lambda_X^{\mathrm{nef}}$, and $\cR_X$ is infinite. 
Hence \eqref{eq:amp-boundary} gives
\begin{equation}\label{eq:round-complement}
\SSph^{\rho(X)-2}\setminus\Lambda_X
=\bigcup_{C\in\cR_X}B_C.
\end{equation}
As observed above, each $B_C$ is a nonempty
round open ball.
By \cref{def:sierpinski-type,lem:root-ball-comparison},
their closures are pairwise disjoint.
Thus the $B_C$ are precisely the complementary
components of $\Lambda_X$, and
$C\mapsto B_C$ is a bijection.
Moreover, $\SSph^{\rho(X)-2}\setminus B_C$ is homeomorphic
to a closed $(\rho(X)-2)$-ball.

By \cref{proper-lim-measure-zero},
$\Lambda_X$ has zero spherical measure.
Since every nonempty open subset of
$\SSph^{\rho(X)-2}$ has positive spherical measure,
the union in \eqref{eq:round-complement}
is dense in $\SSph^{\rho(X)-2}$.

Since distinct $(-2)$-curves have distinct classes
in $\NS(X)$, the set $\cR_X$ is countably infinite.
Enumerate the balls as $B_1,B_2,\ldots$.
For every $\varepsilon>0$, round open balls of
spherical diameter at least $\varepsilon$
have uniformly positive spherical measure.
Since the $B_i$ are disjoint and $\SSph^{\rho(X)-2}$
has finite measure, it follows that
$\diam(B_i)\to0$.
\end{proof}

\begin{question}[Converse to \cref{thm:round-limit-set}]
Let $X$ be a K3 surface %over an algebraically closed field 
with $\rho(X)\geq4$. 
If
\[
 \Lambda_X\subset\partial_\infty\HH_X\cong\SSph^{\rho(X)-2}
\]
is a $(\rho(X)-3)$-dimensional round \sier\ space, must $X$ be
of \sier\ type?
\end{question}

\subsection{Rigidity}\label{subsec:rigidity}

We prove that, among K3 surfaces of \sier\ type, 
the quasi-isometry class of the automorphism group determines the rational N\'eron--Severi lattice together with its integral root set.

%Throughout this subsection, K3 surfaces may be defined over different algebraically closed fields. 
%We use the notation of Section~\ref{sec:preliminaries} and Subsections~\ref{subsec:fundamental}--\ref{subsec:characterization}.
%In particular, $\NS(X)$, $\Delta_X$ is its set of integral roots, and $\cR_X$ is the set of classes of $(-2)$-curves.

We first establish the rational spanning property needed to recover
the rational N\'eron--Severi lattice from the peripheral spheres.

\begin{lemma}\label{lem:roots-span}
Let $X$ be a K3 surface of \sier\ type. Then
\begin{equation}\label{eq:roots-span}
 \NS(X)_\QQ
 =\spanop_{\QQ}\cR_X
 =\spanop_{\QQ}\Delta_X.
\end{equation}
\end{lemma}

\begin{proof}
Put $V:=\spanop_{\QQ}\cR_X$ and regard $V_{\RR}:=V\otimes_{\QQ}\RR$ as a subspace of $\NS(X)_\RR$.
The group $\Gamma_X$ preserves $V$, since automorphisms permute the $(-2)$-curves. 
Choose distinct $C,D\in\cR_X$.
Since $(C,D)\geq3$ by the definition of \sier\ type,
we have $(C+D)^2>0$. 
%\[
% (C+D)^2=-4+2(C,D)>0.
%\]
$C+D$ belongs to $\cC_X$ because it is effective.
Thus
$o:=(C+D)/\sqrt{(C+D)^2}$ is a point of
$\HH_X\cap V_{\RR}$.
Taking accumulation points of $\Gamma_Xo$, we obtain
\begin{equation}\label{eq:limit-in-root-span}
 \Lambda_X\subset
 \{[v]\in\partial_\infty\HH_X\mid v\in V_{\RR}\}.
\end{equation}

For any $C\in\cR_X$, Theorem~\ref{thm:round-limit-set} gives
\[
 \partial B_C
 =\{[v]\in\partial_\infty\HH_X\mid(v,C)=0\}
 \subset\Lambda_X.
\]
Consequently, every nonzero isotropic vector in
$C^\perp\cap\overline{\cC_X}$ belongs to $V_{\RR}$
by \eqref{eq:limit-in-root-span}.
The form on $C^\perp$ has signature
$(1,\rho(X)-2)$, so these isotropic vectors span $C^\perp$
over $\RR$.
Indeed, in Lorentz-orthonormal coordinates $e_0,e_1,\cdots,e_{\rho(X)-2}$ with
$e_0\in\cC_X$, the isotropic vectors $e_0+e_j$ and $e_0-e_j$,
$(1\leq j\leq\rho(X)-2)$, already span $C^\perp$.
It follows that $C^\perp\subset V_{\RR}$.
Since $C\in V_{\RR}$ and $C^2=-2$, we have
$
 \NS(X)_\RR=C^\perp\oplus\RR C=V_{\RR}.
$
Hence $\dim_{\QQ}V=\rho(X)$ and $V=\NS(X)_\QQ$.
The inclusions $\cR_X\subset\Delta_X\subset \NS(X)_\QQ$ give
the second equality in \eqref{eq:roots-span}.
\end{proof}

We now compare the Bowditch boundary with the limit set for K3 surfaces, specifying the metric identification relevant to the rigidity theorem of Bonk--Kleiner--Merenkov (\cref{thm:bkm}).
We also verify the spherical measure condition required to apply this theorem. 

Recall that if $X$ is a K3 surface satisfying \eqref{eq:char}
and $\Aut(X)$ is non-elementary,
then $(\Aut(X),\cP_X)$ is relatively hyperbolic by \cref{thm:k3-relative}.
Here $\cP_X$ is a finite set of representatives of the conjugacy classes of maximal parabolic subgroups of $\Aut(X)$ with respect to its action on $\HH_X$.
%We allow $\cP_X=\varnothing$.
%
There is a natural $\Aut(X)$-equivariant homeomorphism
\[
 \Phi_X:\partial(\Aut(X),\cP_X)
 \xrightarrow{\sim}\Lambda_X
\]
induced by the action of $\Aut(X)$ on $\HH_X$ (\cite{Bow12}).
More precisely, if a sequence $g_n\in\Aut(X)$
converges to $\xi\in\partial(\Aut(X),\cP_X)$
in the compactification of the cusped Cayley graph,
then
$
 \Phi_X(\xi):=\lim_{n\to\infty}g_no
$
for any fixed $o\in\HH_X$.

\begin{proposition}\label{prop:k3-boundary}
Let $X$ be a K3 surface satisfying \eqref{eq:char} and that
$\Aut(X)$ is non-elementary. 
Then the homeomorphism 
$\Phi_X:\partial(\Aut(X),\cP_X)\xrightarrow{\sim}\Lambda_X$ 
is quasisymmetric
for a visual metric on the Bowditch boundary and
the spherical metric on the limit set.
\end{proposition}

\begin{proof}
The assertion is well known to experts,
but we include a proof for completeness.
Let $G:=\Orth^+(\NS(X))$ and $H:=\Gamma_X$.
The group $G$ is a lattice in $\operatorname{Isom}(\HH_X)$, 
hence geometrically finite.
The subgroup $H$ is also geometrically finite by \cite[Theorem~A]{Kik26}.

Let $\cP_G$ be a finite set of representatives of the
conjugacy classes of maximal parabolic subgroups of $G$.
By \cite[Corollary~1.3]{Hru10},
$H$ is relatively quasiconvex in $(G,\cP_G)$.
The peripheral structure $\cP_H$ induced from $\cP_G$
consists of representatives of the $H$-conjugacy classes
of infinite intersections
\[
 H\cap gPg^{-1},
 \qquad g\in G,\quad P\in\cP_G.
\]
These are precisely the maximal parabolic subgroups
of $H$, and $(H,\cP_H)$ is relatively hyperbolic
by \cite[Theorem~9.1]{Hru10}.

Since $G$ and $H$ are finitely generated,
by \cite[Theorem~1.1(2),(3) and Corollary~9.6]{HH21},
the inclusion $H\subset G$ induces a quasisymmetric embedding
%orbit mapで定めるという定義から(roughly H-equivariantではなく)H-equivariantになっている
\[
 \partial(H,\cP_H)\rightarrow
 \partial(G,\cP_G)
 \qs\partial_\infty\HH_X
 %\xrightarrow{\sim_{\mathrm{qs}}}\partial_\infty\HH_X %ここもquasisymmetric
\]
for visual metrics.
More precisely, if $h_n\in H$ converges to
$\xi\in\partial(H,\cP_H)$ in the cusped compactification,
then the boundary map sends $\xi$ to
$\lim_{n\to\infty}h_no\in\partial_\infty\HH_X$, where $o\in\HH_X$ is fixed. 
Its image is therefore $\Lambda_H=\Lambda_X$.
As explained in \cref{subsec:round-limit-sets}, we have 
\[
\partial_\infty\HH_X\qs\SSph^{\rho(X)-2}
\]
with respect to a visual metric on $\partial_\infty\HH_X$ and the spherical metric on $\SSph^{\rho(X)-2}$. 
%Moreover, a visual metric on $\partial_\infty\HH_X$ is quasisymmetric equivalent to the spherical metric by \eqref{eq:spherical-visual}.

Finally, the homomorphism $\Aut(X)\to H$ has finite
kernel and identifies $\cP_X$ with the inverse images
of $\cP_H$.
It is therefore a quasi-isometry of pairs and induces
a quasisymmetry of their Bowditch boundaries
$\partial(\Aut(X),\cP_X)\qs\partial(H,\cP_H)$
by \cref{thm:relative-qi}.
The resulting identification with $\Lambda_X$
agrees with $\Phi_X$.
\end{proof}

\begin{theorem}[Quasi-isometry rigidity]%[Theorem D]
\label{thm:rigidity-sierpinski}
Let $X_1$ and $X_2$ be K3 surfaces of \sier\ type.
The following conditions are equivalent:
\begin{enumerate}
\item\label{item:rigidity-qi}
$\Aut(X_1)\qi\Aut(X_2)$;
\item\label{item:rigidity-lattice}
there is a root-preserving rational isometry
$\phi:\NS(X_1)_\QQ\xrightarrow{\sim}\NS(X_2)_\QQ$;
\item\label{item:rigidity-mobius}
$\rho(X_1)=\rho(X_2)$, and their limit sets are \mobius\ equivalent in the corresponding boundary sphere.
\end{enumerate}

\begin{comment}
When these conditions hold, $\phi$ may be chosen so that
\begin{equation}\label{eq:rigidity-ample}
 \phi(\cR_{X_1})=\cR_{X_2},\qquad
 \phi_{\RR}(\Amp(X_1))=\Amp(X_2).
\end{equation}
For such a choice, $\phi\Gamma_{X_1}\phi^{-1}$ and
$\Gamma_{X_2}$ are commensurable subgroups of
$\Orth^+(\NS(X_2)_\QQ)$.
%Subsection 2.2のPropとダブるのでカット
\end{comment}
\end{theorem}

\begin{proof}
Write $L_i:=\NS(X_i)$, $\rho_i:=\rho(X_i)$,
$\Lambda_i:=\Lambda_{X_i}$, and $\cP_i:=\cP_{X_i}$.
Suppose first that \textup{\ref{item:rigidity-qi}} holds,
and let $f:\Aut(X_1)\to\Aut(X_2)$ be a quasi-isometry.
We show that $f$ induces a quasisymmetry
\begin{equation}\label{eq:limit-quasisymmetry}
 f_\Lambda:\Lambda_1\xrightarrow{\sim}\Lambda_2
\end{equation}
for the spherical metrics.

By \cref{thm:fundamental}~\ref{item:hyperbolic}
and \cref{thm:k3-relative},
$\cP_i=\varnothing$ if and only if $\Aut(X_i)$
is hyperbolic.
Since hyperbolicity is invariant under quasi-isometry,
either both peripheral families are empty or neither is.
In the latter case, every member of $\cP_i$ is virtually
$\ZZ^{\rho_i-2}$ by
\cref{thm:fundamental}~\ref{item:fibrations}
and \cref{thm:k3-relative}.
Since $\rho_i\geq4$, these groups have rank at least two.
Thus $f$ coarsely preserves the peripheral families
by \cref{thm:nrh-peripherals}.
In the former case, this is automatic.
Hence $f$ is a quasi-isometry of pairs, and
\cref{thm:relative-qi,prop:k3-boundary}
give the quasisymmetry
\[
 f_\Lambda:=\Phi_{X_2}\circ\partial f\circ\Phi_{X_1}^{-1}
\]
in \eqref{eq:limit-quasisymmetry}.

By \cref{thm:round-limit-set},
$\Lambda_i$ is a Schottky set of 
topological dimension $\rho_i-3$.
Thus \eqref{eq:limit-quasisymmetry} implies
$\rho_1=\rho_2$.
Put $m:=\rho_1-2\geq2$ and identify both ideal
boundaries with $\SSph^m$.
%By \cref{thm:round-limit-set,proper-lim-measure-zero}, the sets $\Lambda_i$ are Schottky sets of spherical measure zero.
By \cref{proper-lim-measure-zero}, $\Lambda_i$ is of spherical measure zero.
Therefore, the Bonk--Kleiner--Merenkov rigidity theorem
(\cref{thm:bkm}) implies that $f_\Lambda$ is the
restriction of a \mobius\ transformation of $\SSph^m$.
This proves \textup{\ref{item:rigidity-mobius}}.

Next, suppose that \textup{\ref{item:rigidity-mobius}}
holds.
Lift the given \mobius\ transformation
$M:\SSph^m\to\SSph^m$ with $M(\Lambda_1)=\Lambda_2$
to a hyperbolic isometry
$F:\HH_{X_1}\to\HH_{X_2}$ with $\partial_\infty F=M$
and a linear isometry
\[
 T:L_{1,\RR}\xrightarrow{\sim}L_{2,\RR},
 ~\text{ with }~T(\cC_{X_1})=\cC_{X_2}
 ~\text{ and }~T|_{\HH_{X_1}}=F.
\]
Since $\partial_\infty F$ maps $\Lambda_1$ onto $\Lambda_2$,
it bijectively maps their complementary components.
By \cref{thm:round-limit-set}, for each
$C\in\cR_{X_1}$ there is a unique
$D\in\cR_{X_2}$ such that $\partial_\infty F(B_C)=B_D$.
Since $F$ is a hyperbolic isometry,
we have $T(C^\perp)=D^\perp$.
It follows that $T(C)=aD$ for some $a>0$.
Since $C^2=D^2=-2$ and $T$ is an isometry,
we have $a=1$.
Applying the same argument to $T^{-1}$ gives
\begin{equation}\label{eq:curve-bijection}
 T(\cR_{X_1})=\cR_{X_2}.
\end{equation}
%In particular, \eqref{eq:ample} implies that
%\[
% T(\Amp(X_1))=\Amp(X_2).
%\]
By \cref{lem:roots-span}, %and \eqref{eq:curve-bijection},
we have
\begin{equation}\label{eq:rational-descent}
 T(L_{1,\QQ})
 =T(\spanop_{\QQ}\cR_{X_1})
 =\spanop_{\QQ}\cR_{X_2}
 =L_{2,\QQ}.
\end{equation}
\begin{comment}
For any K3 surface $Y$, the Weyl group is generated by the
reflections in the walls of the ample chamber, and
\begin{equation}\label{eq:roots-weyl-orbits}
 \Delta_Y=\{\varepsilon w(C)\mid
      \varepsilon\in\{1,-1\},\ w\in W_Y,\ C\in\cR_Y\}.
\end{equation}
Indeed, local finiteness of the root arrangement inside the
positive cone makes every root hyperplane a wall of some
chamber. An element of $W_Y$ takes this chamber to the ample
chamber and its normal root to a $(-2)$-curve class up to sign;
see \cite[Chapter~8, \S2]{Huy16}.-
\end{comment}
Moreover, the Weyl chamber decomposition of positive cones recalled in Subsection~2.1 implies that
every integral root is a Weyl group translate of a $(-2)$-curve class up to sign.
Since we have
$
 T s_C T^{-1}=s_{T(C)},
$
equation \eqref{eq:curve-bijection} implies that
$T W_{X_1}T^{-1}=W_{X_2}$ and hence
$T(\Delta_{X_1})=\Delta_{X_2}$.
Thus $\phi:=T|_{L_{1,\QQ}}$ is a root-preserving
rational isometry.  %satisfying \eqref{eq:rigidity-ample}.
This proves \textup{\ref{item:rigidity-lattice}}.

Finally, the implication
\textup{\ref{item:rigidity-lattice}} $\Rightarrow$
\textup{\ref{item:rigidity-qi}} follows from
\cref{prop:root-qi}, since $\rho_i\in\{4,5\}$.
%The commensurability assertion follows from the chamber-stabilizer argument in the proof of \cref{prop:root-qi}, using \cref{thm:cone-comparison}.
\end{proof}

% !TEX root = ../main.tex
% !TEX root = ../main.tex
% Integrated from section4_1_hyperbolic_qi_v2.tex; full change log: ../CHANGELOG-ja.md.
\section{Infinite quasi-isometry classes for Picard rank 4}\label{sec:rank-four}
\subsection{Quasi-isometry classes of hyperbolic automorphism groups}
\label{subsec:hyperbolic-qi-classes}

Hyperbolic groups form a basic class of negatively curved groups in geometric group theory. 
They include free groups of finite rank and fundamental groups of closed hyperbolic manifolds. 
Hyperbolicity is invariant under quasi-isometry, and Gromov boundaries of hyperbolic groups retain information on their large-scale geometry.
In particular, a quasi-isometry induces a homeomorphism of the
boundaries, and a quasisymmetry when they are equipped with visual metrics. 
Thus boundary geometry provides a natural approach to the quasi-isometry classification. %; see Section~\ref{sec:preliminaries} and \cite[Chapter~III.H]{BH99}.
See \cref{subsec:relhyp} for details. 
In this subsection, we consider this problem for hyperbolic
automorphism groups of K3 surfaces.

\begin{definition}\label{def:Q-rho}
For $\rho\in\{1,\ldots,20,22\}$, we set
\[
 \cQ_\rho:=
 \left\{
 [\Aut(X)]_{\mathrm{qi}}\ \middle|\
 \begin{array}{l}
 X/k\text{ is a K3 surface},\\
 \rho(X)=\rho,\\
 \Aut(X)\text{ is non-elementary hyperbolic.}
 \end{array}
 \right\},
 \qquad Q_\rho:=\#\cQ_\rho.
\]
Here $k$ ranges over algebraically closed fields of arbitrary
characteristic, and $[G]_{\mathrm{qi}}$ denotes the quasi-isometry
class of a finitely generated group $G$ with a word metric.
\end{definition}

The non-elementary condition excludes finite and virtually infinite
cyclic groups. For Picard ranks at most twenty, the collection of
classes can be computed over $\CC$.

\begin{proposition}\label{prop:complex-reduction}
Let $X$ be a K3 surface over an algebraically closed field $k$ with
$1\leq\rho(X)\leq20$. Then there exists a complex 
%projective 
K3 surface $Y$ such that
\[
 \NS(X)\cong\NS(Y),\qquad \Aut(X)\qi\Aut(Y).
\]
Consequently, for $1\leq\rho\leq20$, the set $\cQ_\rho$ is
unchanged if Definition~\ref{def:Q-rho} is restricted to $k=\CC$.
\end{proposition}

\begin{proof}
Suppose first that $\charac k=0$.
Choose a subfield $k_0\subset k$ finitely generated over $\QQ$ such that $X$ is defined over $k_0$. %descend X to k_0
Since $k_0$ is finitely generated over $\QQ$,
it admits an embedding into $\CC$.
By the invariance of N\'eron--Severi lattices under extensions of algebraically closed fields \cite[Proposition~3.1]{MP12},
the resulting complex K3 surface $Y$ satisfies $\NS(X)\cong\NS(Y)$.

Suppose next that $\charac k=p>0$.
Since $\rho(X)\leq20$, $X$ has finite height. 
%the Tate conjecture for K3 surfaces implies that $X$ has finite height \cite{MP15,KM16,MP20}+Ito-Ito-Koshikawa(p=2)
%Yuya Matsumoto, On μ_n-actions on K3 surfaces in positive characteristic⁠, Nagoya Math. J. 249 (2023), 11–49.のLemma 8.3もよくまとまってる
By \cite[Corollary~4.2 and Section~6]{LM18},
$X$ admits a Picard-preserving lifting to characteristic zero, including when $p=2$. 
Applying the characteristic-zero case to its geometric generic fiber yields the required $Y$.

%By \cite[Corollary~4.2]{LM18}, 
%there is a smooth projective lifting over the Witt ring $W(k)$ whose geometric generic fiber $X_{\bar\eta}$ satisfies
%\[
% \NS(X_{\bar\eta})\xrightarrow{\ \sim\ }\NS(X).
%\]
%Indeed, all classes on the special fiber lift, and specialization is injective. 
%The intersection pairing is preserved. This lifting statement also holds for $p=2$, as explicitly noted in \cite[Remark~6.3]{LM18}. Applying the characteristic-zero argument to $X_{\bar\eta}$ gives the required complex surface $Y$.

In both cases, the above isomorphisms preserve the intersection pairing and hence the root sets.
Since both K3 surfaces satisfy \textup{(2.3)}, Proposition~\ref*{prop:root-qi} gives $\Aut(X)\qi\Aut(Y)$.
The last assertion follows from the quasi-isometry invariance of hyperbolicity and non-elementarity.
\end{proof}

\begin{remark}
For $1\leq\rho\leq20$, the set $\cQ_\rho$ is therefore at most
countable: there are only countably many integral lattices, and
isometric N\'eron--Severi lattices give quasi-isometric automorphism
groups by Proposition~\ref{prop:root-qi}. Together with
$Q_{22}=0$ in \cref{prop:Q-22} below, this shows
that $Q_\rho=\infty$, when it occurs, means that there are countably infinitely many classes.
\end{remark}

\begin{comment}
\begin{remark}\label{rem:base-field-Q}
The algebraic closure in Proposition~\ref{prop:complex-reduction}
is essential to its formulation. If $X$ is defined over a field $k$
which is not algebraically closed, the proposition applies to
$\Aut(X_{\bar k})$ and the geometric Picard rank
$\rho(X_{\bar k})$. It does not assert the same conclusion for
$\Aut_k(X)$, since automorphisms need not descend to $k$.
%
Also, allowing all algebraically closed fields to vary does not
assert that every class occurs over each fixed field of positive
characteristic. 
%
The supersingular case $\rho=22$ cannot occur over
$\CC$ and is treated separately in Proposition~\ref{prop:Q-22}.
\end{remark}
\end{comment}

We recall two classical characterizations by the topology of the
Gromov boundary. Here a surface group means the fundamental group of a closed orientable (real) surface of genus at least two, and $F_r$ denotes the free group of rank $r$.

\begin{theorem}[See {\cite[p.648]{KK00}}]\label{thm:boundary-free-surface}
Let $G$ be a non-elementary hyperbolic group.
\begin{enumerate}
\item $\dim\partial G=0$ if and only if $G$ is virtually free.
In this case, $\partial G$ is a Cantor set and $G\qi F_2$.
\item $\partial G$ is homeomorphic to $\SSph^1$ if and only if
$G$ is virtually a surface group. In this case, $G\qi\HH^2$.
\end{enumerate}
\end{theorem}

The first characterization follows from the classical theory of
ends and splittings of hyperbolic groups. The second follows from
the convergence-group theorems of Tukia, Gabai, and
Casson--Jungreis. 
%We refer to \cite[p.~648]{KK00} for these statements and the original references. 
For the quasi-isometry assertions, all nonabelian free groups of finite rank are commensurable, 
and all closed hyperbolic surface groups act properly and cocompactly on $\HH^2$, 
so the conclusions follow from the \v Svarc--Milnor lemma.
In particular, these are two distinct quasi-isometry classes, since a Cantor set is not homeomorphic to $\SSph^1$.

\begin{proposition}\label{prop:Q-low-rank}
Let $X$ be a K3 surface. %over an algebraically closed field.
\begin{enumerate}
\item If $\rho(X)=1$, then $\Aut(X)$ is finite.
\item If $\rho(X)=2$, then $\Aut(X)$ is virtually cyclic.
\item If $\rho(X)=3$, then $\Aut(X)$ is hyperbolic.
If it is non-elementary, then
$
\Aut(X)\qi F_2
$
or 
$
\Aut(X)\qi \HH^2,
$
and both cases occur.
\end{enumerate}
Consequently, we have
\[
Q_1=Q_2=0,\qquad Q_3=2.
\]
\end{proposition}

\begin{proof}
Since $\rho(X)\leq3$, it suffices to work over $\CC$
by Proposition~\ref{prop:complex-reduction}. 
Recall that, by Theorem~\ref{thm:cone-comparison},
the representation $\varrho_X:\Aut(X)\twoheadrightarrow\Gamma_X$ has finite kernel,
and both are finitely generated.
If $\rho(X)=1$, then $\Gamma_X$ is trivial, so
$\Aut(X)$ is finite.
If $\rho(X)=2$, then $\Gamma_X$ is a discrete
subgroup of $\Isom(\HH^1)$ and hence virtually
cyclic. Thus so is $\Aut(X)$.

Suppose that $\rho(X)=3$.
Then $\Gamma_X$ is a finitely generated discrete
subgroup of $\Isom(\HH^2)$.
If $\Gamma_X$ is elementary, it is virtually cyclic.
Otherwise, by Selberg's lemma, it has a torsion-free
subgroup $G$ of finite index contained in
$\Isom^+(\HH^2)$.
By the standard presentation of finitely generated
Fuchsian groups 
(see, e.g., \cite{Kat92, BS83}),
%(see, e.g., \cite[Sections~4.3 and~4.6]{Kat92} and \cite[p.~414]{BS83}),
it is well known that every non-elementary, torsion-free, finitely generated Fuchsian group is either a nonabelian free group of finite rank or a closed hyperbolic surface group. 
%By \cite[p.~25, proof of Theorem~1.1]{BCR16}, $G$ is either a finitely generated free group or a closed hyperbolic surface group.
Since $\varrho_X$ has finite kernel,
$\Aut(X)\qi\Gamma_X\qi G$.
Thus $\Aut(X)$ is hyperbolic, and if non-elementary,
it is quasi-isometric to either $F_2$ or $\HH^2$.
Both classes occur by the examples below, and
they are distinct by
Theorem~\ref{thm:boundary-free-surface}.
\end{proof}

\begin{example}[Virtually free groups]\label{ex:rank-three-free}
Let $X$ be a very general smooth hypersurface of multidegree
$(2,2,2)$ in $\PP^1\times\PP^1\times\PP^1$.
Then $X$ is a Wehler K3 surface with $\rho(X)=3$ and
\[
 \NS(X)\cong
 \begin{pmatrix}
 0&2&2\\2&0&2\\2&2&0
 \end{pmatrix}
 \cong U(2)\oplus\langle-4\rangle.
\]
The three projections to $\PP^1\times\PP^1$ are double covers,
and their covering involutions generate
\[
 \Aut(X)\cong C_2*C_2*C_2,
\]
where $C_2$ is the cyclic group of order 2; see
\cite[Introduction and Example~1.5]{HL25}.
The kernel of the homomorphism $C_2*C_2*C_2\to C_2$ which is
nontrivial on each factor is a free group of rank 2 and index 2.
Thus $\Aut(X)\qi F_2$ and $\partial\Aut(X)$ is a Cantor set. 

These examples predate the work of Hashimoto--Lee. 
For explicit Wehler surfaces over $\QQ$ with geometric Picard rank 3, see \cite[Thorem 2.1]{BvL07}. 
Hashimoto--Lee construct K3 surfaces with
\[
 \NS(X_n)\cong U(n)\oplus\langle-2n\rangle\qquad(n\geq2)
\]
and determine their automorphism groups as congruence subgroups of
$\operatorname{PGL}_2(\ZZ)$; see \cite[Proposition~1.4]{HL25}.
In particular, they obtain free automorphism groups of unbounded
finite rank \cite[Corollary~1.8]{HL25}. All these nonabelian free
groups represent the same element of $\cQ_3$.
\end{example}

\begin{comment}
\begin{example}%\label{ex:rank-three-surface}
Consider the lattice
\[
 L:=\langle4\rangle\oplus\langle-8\rangle\oplus\langle-12\rangle.
\]
It contains no roots, since every square is divisible by four.
It also contains no nonzero isotropic vector. Indeed, an integral
solution of
\[
 x^2-2y^2-3z^2=0
\]
satisfies $x\equiv y\equiv0\pmod3$, and then $z\equiv0\pmod3$.
Infinite descent excludes a nonzero solution.
The lattice $L$ embeds primitively into $U^{\oplus3}$: in standard
bases $e_i,f_i$ of the three copies of $U$, take
$e_1+2f_1$, $e_2-4f_2$, and $e_3-6f_3$.
By the surjectivity of the period map, there exists a complex
projective K3 surface $X$ with $\NS(X)\cong L$.

Since $L$ has no roots, the ample cone is the positive cone.
Thus $\Gamma_X$ has finite index in $\Orth^+(L)$.
The anisotropy of $L$ implies that this arithmetic subgroup of
$\Isom(\HH^2)$ is cocompact. Hence $\Aut(X)$ is virtually a closed surface group and
\[
 \Aut(X)\qi\HH^2,\qquad \partial\Aut(X)\cong\SSph^1.
\]
This is the case $k=1$ of \cite[Example~4.21]{Kik26}; see also
\cite[Corollary~4.13\textup{(ii)}]{Kik26}.
\end{example}
\end{comment}

\begin{example}[Virtually surface groups]\label{ex:rank-three-surface}
For each positive integer $k\equiv1\pmod3$, 
let $X_k$ be a complex K3 surface with
\begin{equation}\label{eq:rank-four-2-manifold-lattice}
 \NS(X_k)\cong 
 \langle4\rangle\oplus\langle-8\rangle
 \oplus\langle-12k\rangle.
\end{equation}
By \cite[Corollary~4.13\textup{(ii)} and Example~4.21]{Kik26}, $\Aut(X_k)$ is virtually a closed surface group. 
Consequently,
\[
 \Aut(X_k)\qi\HH^2,
 \qquad
 \partial\Aut(X_k)\cong\SSph^1.
\]
\end{example}

\begin{remark}\label{rem:fuchsian-boundaries}
The Gromov boundary $\partial\Aut(X)$ need not be
homeomorphic to the limit set $\Lambda_X\subset\partial_\infty\HH_X$. %of its natural representation. 
In Example~\ref{ex:rank-three-free}, the lattice $U(2)\oplus\langle-4\rangle$ has no roots and an isotropic vector. 
Thus $\Gamma_X$ is a nonuniform lattice in $\Isom(\HH^2)$, and $\Lambda_X=\SSph^1$, whereas $\partial\Aut(X)$ is a Cantor set. 
%This is why the proof of Proposition~\ref{prop:Q-low-rank} uses the classification of finitely generated Fuchsian groups, rather than identifying the two boundaries.
\end{remark}

For higher Picard ranks, a finiteness theorem is provided by
Fujiwara--Oguiso--Yu.

\begin{theorem}[Fujiwara--Oguiso--Yu,~{\cite[Theorem~1.8 and Table~1]{FOY26}}]\label{thm:FOY-finite-lattices}
The set of isometry classes of N\'eron--Severi lattices of complex K3 surfaces $X$ such that $\rho(X)\geq6$ and $\Aut(X)$
is non-elementary hyperbolic consists of exactly $166$ lattices.
\end{theorem}

%The number $166$ counts lattices over all Picard ranks at least six, not at each fixed rank.

\begin{corollary}\label{cor:Q-high-rank}
%For every $6\leq\rho\leq20$, we have $Q_\rho\leq166$.
%More precisely,
We have 
\[
 \sum_{\rho=6}^{20}Q_\rho\leq166.
\]
\end{corollary}

\begin{proof}
By Proposition~\ref{prop:complex-reduction}, it suffices to work
over $\CC$.
For each $\rho$, let $\mathcal L_\rho$ denote the subset of the
$166$ lattices in Theorem~\ref{thm:FOY-finite-lattices} having
rank $\rho$. Proposition~\ref{prop:root-qi} gives a well-defined
surjection
\[
 \mathcal L_\rho\longrightarrow\cQ_\rho,\qquad
 [\NS(X)]\longmapsto[\Aut(X)]_{\mathrm{qi}}.
\]
Hence $Q_\rho\leq\#\mathcal L_\rho$, and summing proves the
assertion.
\end{proof}

%Distinct lattices may give the same quasi-isometry class, so the lattice count does not determine the numbers $Q_\rho$.

\begin{proposition}\label{prop:Q-22}
We have $Q_{22}=0$. 
%Automorphism groups of K3 surfaces of Picard rank 22 are not hyperbolic. Therefore, we have $Q_{22}=0$. 
\end{proposition}

\begin{pf}
%Suppose that $\rho(X)=22$. Then 
$X$ is supersingular in positive characteristic.
By \cref{lem:supersingular-mw},
$\Aut(X)$ contains $\ZZ^2$ and hence is not hyperbolic. 
%Thus $Q_{22}=0$.
\end{pf}

The preceding results leave the finiteness question for $Q_4$ and
$Q_5$. 
We turn to the case of Picard rank 4 in the next subsection.
% !TEX root = ../main.tex
% Integrated from section4_2_picard_rank_four(1).tex; full change log: ../CHANGELOG-ja.md.
\subsection{Case of Picard rank 4}\label{subsec:picard-rank-four}

We now show that the finiteness result in higher Picard rank has no
analogue in Picard rank 4. 
We use the notation $\cQ_\rho$ and $Q_\rho$ of Definition~\ref{def:Q-rho}. 
All K3 surfaces constructed in this subsection are over $\CC$: %complex K3 surfaces. 
recall that every even lattice of signature $(1,3)$ occurs as their N\'eron--Severi lattice. 
%, by the primitive embedding theorem and the surjectivity of the period map \cite[Corollary~14.3.1 and Theorem~6.3.1]{Huy16}. 

\begin{theorem}%[Theorem A]
\label{thm:Q-four-infinite}
We have $Q_4=\infty$. More precisely, there exist infinitely many
complex K3 surfaces of \sier\ type and Picard rank 4 whose
automorphism groups are non-elementary hyperbolic and pairwise
non-quasi-isometric.
\end{theorem}

We first extend the lattice in \cref{ex:Sierpinski-lattice} (i) to a
one-parameter family. Set
\[
 \mathcal T:=\{t\in\ZZ_{>0}\mid ~t\equiv3\pmod{12}\}.
\]
For $t\in\mathcal T$, let $L_t$ be the lattice with Gram matrix
\begin{equation}\label{eq:rank-four-family}
 G_t:=
 \begin{pmatrix}
 2&t&t&t\\
 t&-2&4&4\\
 t&4&-2&4\\
 t&4&4&-2
 \end{pmatrix}.
\end{equation}

\begin{lemma}\label{lem:rank-four-family}
For every $t\in\mathcal T$, the lattice $L_t$ satisfies the followings. 
\begin{enumerate}
\item
$L_t$ is even hyperbolic with determinant $-108(t^2-4)$. 
%has signature $(1,3)$ and determinant $-108(t^2-4)$. 
\item
$L_t$ is anisotropic. 
\item 
$L_t$ contains roots, and any two roots $r\ne\pm s$ satisfy
$|(r,s)|\geq3$.
\end{enumerate}
\end{lemma}

\begin{proof}
(i)~The diagonal entries are even, so $L_t$ is even.
The lower-right block of $G_t$ has eigenvalues $6,-6,-6$.
Its Schur complement is $2-t^2/2<0$. This gives signature $(1,3)$ and
\[
 \det G_t=216(2-t^2/2)=-108(t^2-4).
\]

\noindent
(ii)~For $v=(x,y_1,y_2,y_3)\in L_t$, put $S:=y_1+y_2+y_3$ and $T:=y_1^2+y_2^2+y_3^2$. 
Then
\begin{equation}\label{eq:rank-four-family-form}
 q_t(v):=\tfrac12(v,v)=x^2+txS+2S^2-3T.
\end{equation}
Suppose that $q_t(v)=0$. 
Since $t$ is odd, %reduction modulo two gives 
we have $x+xS+S\equiv0\pmod2$, which forces both $x$ and $S$ to be even.
If some $y_i$ is odd, exactly two are odd, and $T\equiv2\pmod4$.
The other three terms in \eqref{eq:rank-four-family-form} are
divisible by 4, so $q_t(v)\equiv2\pmod4$, a contradiction.
Thus all four coordinates are even. Infinite descent excludes a
nonzero integral isotropic vector. %, and hence a nonzero rational one.

\noindent
(iii)~
The last three basis vectors of $L_t$ are roots. 
For a root $r$, the equation $q_t(r)\equiv-1\pmod3$ %modulo three 
gives
$x_r\equiv0$ and $S_r\not\equiv0\pmod3$. 
%Polarizing
By \eqref{eq:rank-four-family-form}, we obtain, for any two roots,
\[
 (r,s)\equiv S_rS_s\not\equiv0\pmod3.
\]
In particular, their pairing is not zero. 
If $(r,s)=1$, then $u=r-s$ satisfies $q_t(u)=-3$ and $x_u\equiv S_u\equiv0\pmod3$.
%Reduction modulo nine gives $T_u\equiv1\pmod3$. 
By \eqref{eq:rank-four-family-form}, we have $-3\equiv q_t(u)\equiv -3T_u\pmod9$, hence $T_u\equiv1\pmod3$. 
Among three coordinates, this means that exactly one $y_i$ is nonzero modulo 3, contradicting $S_u\equiv0\pmod3$. Replacing $s$ by $-s$
also excludes $(r,s)=-1$. Finally, $|(r,s)|=2$ for $r\ne\pm s$
would make one of $r+s,r-s$ a nonzero isotropic vector.
\end{proof}

\begin{cor}\label{cor:rank-four-family}
For every $t\in\mathcal T$, 
a complex K3 surface $X_t$ with $\NS(X_t)\cong L_t$ is of \sier\ type and has a non-elementary hyperbolic automorphism group.
\end{cor}
\begin{pf}
The assertion follows from \cref{lem:rank-four-family}, \cref{lem:numerical} and \cref{thm:fundamental}~(iv)(vi).
\end{pf}

The next elementary observation supplies the required infinite set
of parameters. %, without a squarefree-value assertion for a polynomial.

\begin{lemma}\label{lem:distinct-square-classes}
There is an infinite subset $\mathcal S\subset\mathcal T$ such that
\[
 \frac{t^2-4}{u^2-4}\notin(\QQ^\times)^2
 \qquad\text{for all distinct }t,u\in\mathcal S.
\]
\end{lemma}

\begin{proof}
We construct an increasing sequence $t_1,t_2,\ldots$ in $\mathcal T$.
Take $t_1=3$. 
Given $t_1,\ldots,t_n$ satisfying the property, 
we then choose a prime $p>3$ that does not divide $\prod_{i=1}^n(t_i^2-4)$. By the Chinese remainder theorem, there is an integer $t_{n+1}>t_n$ satisfying
\[
 t_{n+1}\equiv3\pmod{12},\qquad
 t_{n+1}\equiv2+p\pmod{p^2}.
\]
Writing $v_p$ for the $p$-adic valuation, we have
\[
 v_p(t_{n+1}-2)=1,\qquad v_p(t_{n+1}+2)=0,
 \qquad v_p(t_i^2-4)=0\quad(1\leq i\leq n).
\]
Thus $v_p((t_{n+1}^2-4)/(t_i^2-4))=1$, whereas the valuation of a
rational square is even. The set $\mathcal S=\{t_n\mid n\geq1\}$
has the required property.
\end{proof}

\begin{proof}[Proof of Theorem~\ref{thm:Q-four-infinite}]
For each $t\in\mathcal S$, choose a complex K3 surface $X_t$ with $\NS(X_t)\cong L_t$.
Corollary~\ref{cor:rank-four-family} shows that every $X_t$ is of
\sier\ type and that $\Aut(X_t)$ is non-elementary hyperbolic.
Suppose that $\Aut(X_t)\qi\Aut(X_u)$ for distinct $t,u\in\mathcal S$.
Theorem~\ref{thm:rigidity-sierpinski} gives a
root-preserving rational isometry
\[
 L_t\otimes\QQ\xrightarrow{\ \sim\ }L_u\otimes\QQ.
\]
In particular, the determinant ratio must be a rational square:
indeed, $G_t=A^{\mathsf T}G_uA$ for some $A\in\mathrm{GL}_4(\QQ)$
implies
\[
 \frac{t^2-4}{u^2-4}
 =\frac{\det G_t}{\det G_u}=(\det A)^2\in(\QQ^\times)^2.
\]
This contradicts Lemma~\ref{lem:distinct-square-classes}.
We obtain infinitely many distinct elements of $\cQ_4$.
%As observed in Subsection~\ref{subsec:hyperbolic-qi-classes}, this set is at most countable, so it is countably infinite.
\end{proof}

By Theorem~\ref{thm:fundamental}~(vi) and Theorem~\ref{thm:round-limit-set},
$\partial\Aut(X_t)\cong \Lambda_{X_t}$ is a round \sier\ carpet for every $t\in\mathcal T$.
Thus even within this fixed boundary homeomorphism type,
infinitely many quasi-isometry classes occur.
We finish with two further distinct elements of $\cQ_4$,
represented by $F_2$ and $\HH^3$.

\begin{comment}
In passing between $\Aut(X)$ and $\Gamma_X$, we use the finite kernel of the action on $\NS(X)$.
Moreover, $\Aut(X)$ acts faithfully on $H^2(X,\ZZ)$, so Selberg's
lemma supplies a torsion-free subgroup of finite index. Its action
on $\NS(X)$ is faithful as well.
\end{comment}

\begin{example}[Virtually free groups]\label{ex:rank-four-free}
Baragar's example \cite{Bar17} is a complex K3 surface
$X$ with
\begin{equation}\label{eq:rank-four-free-lattice}
 \NS(X)\cong 
 \begin{pmatrix}
 -2&2&2&4\\
 2&-2&2&4\\
 2&2&-2&0\\
 4&4&0&0
 \end{pmatrix}.
\end{equation}
Its ample cone is described by the Apollonian circle packing.
Takatsu's computation for this example (\cite[Example~5.41]{Tak25}) gives $\vcd\Aut(X)=1$. 
Hence $\Aut(X)$ is virtually free by the Stallings--Swan theorem \cite{Swa69}. 
%この時点ではまだvirtually infinite cyclicの可能性がある
%Hence a torsion-free subgroup of finite index has cohomological dimension one and is free by the Stallings--Swan theorem \cite{Swa69}. 
%The corresponding action on $\HH^3$ is non-elementary.
Since $\Aut(X)$ is non-elementary, 
%because limit set is infinite.
it has a finite-index free subgroup of rank at least two.
Consequently,
\[
 \Aut(X)\qi F_2,
 \qquad \partial\Aut(X)\text{ is a Cantor set}.
\]
Here the Gromov boundary is distinct from the limit set %of the natural action on $\HH^3$
, which has parabolic elements (cf. \cref{rem:fuchsian-boundaries}).
\end{example}

\begin{example}[Virtually hyperbolic $3$-manifold groups]
\label{ex:rank-four-three-manifold}
For each positive integer $k\equiv1\pmod3$, 
let $X_k$ be a complex K3 surface with
\begin{equation}\label{eq:rank-four-3-manifold-lattice}
 \NS(X_k)\cong 
 \langle4\rangle\oplus\langle-8\rangle
 \oplus\langle-12\rangle\oplus\langle-12k\rangle.
\end{equation}
By \cite[Corollary~4.13\textup{(ii)} and Example~4.21]{Kik26},
$\Aut(X_k)$ is virtually the fundamental group of a closed hyperbolic $3$-manifold.
Consequently,
\[
 \Aut(X_k)\qi\HH^3,
 \qquad \partial\Aut(X_k)\cong\SSph^2.
\]
\end{example}

\begin{comment}
The three examples have pairwise non-homeomorphic Gromov boundaries,
and hence represent three distinct elements of $\cQ_4$.
Their boundaries also distinguish them from all the groups with
\sier\ carpet boundary in Theorem~\ref{thm:Q-four-infinite}.
\end{comment}

\begin{remark}
It remains open whether $Q_5$ is finite or infinite.
By \cref{thm:fundamental}~(vi), the automorphism group of a
K3 surface of \sier\ type with Picard rank 5 is not
hyperbolic. 
Thus K3 surfaces of \sier\ type do not contribute to $\cQ_5$. 
%, and our approach in this paper provides no information about $Q_5$.
%Thus K3 surfaces of \sier\ type do not contribute to $\cQ_5$, and our approach in this paper does not address whether $Q_5$ is finite or infinite.
\end{remark}

\bibliographystyle{k3alpha}
\bibliography{math}
\end{document}